\documentclass[12pt,english]{article}
\usepackage[latin9]{inputenc}
\usepackage[letterpaper]{geometry}
\usepackage{color}
\usepackage{babel}
\usepackage{amsmath}
\usepackage{amssymb}
\usepackage{setspace}
\usepackage{esint}
\usepackage[authoryear]{natbib}
\usepackage[unicode=true]
 {hyperref}

\makeatletter
\numberwithin{equation}{section}

\usepackage{sectsty}
\usepackage{pslatex}
\usepackage{babel}
\usepackage{chngcntr}
\usepackage{apptools}
\usepackage{bm}
\usepackage{amsmath}
\AtAppendix{\counterwithin{lemma}{section}}
\AtAppendix{\counterwithin{proposition}{section}}

\sectionfont{\sffamily\large}
\subsectionfont{\sffamily\normalsize}
\subsubsectionfont{\sffamily\small}

\newcounter{eqnum}[section] 
\newtheorem{proposition}{Proposition}
\newtheorem{theorem}{Theorem}
\newtheorem{example}{Example}
\newtheorem{lemma}{Lemma} 

\newtheorem{condition}{Condition}

\newcommand{\eproof}{\mbox{}\hfill{\rule{8pt}{8pt}}}

\makeatother

\begin{document}
\title{The Absence of Attrition in \\a War of Attrition under Complete Information\thanks{We thank the editor, Nicolas Vieille, the advisory editor, and two
anonymous referees, Hector Chade, Johannes Hörner, Mike Powell, and
Luis Rayo, as well as participants at several seminars and conferences
for helpful comments and suggestions. A previous version of this paper
was titled ``Equilibrium Selection in the War of Attrition under
Complete Information''.}}
\author{George Georgiadis, Youngsoo Kim, and H. Dharma Kwon\thanks{G. Georgiadis: Kellogg School of Management, Northwestern University;
\protect\href{mailto:}{g-georgiadis@kellogg.northwestern.edu}.; \protect \\
Y. Kim: Culverhouse College of Business, University of Alabama; \protect\href{mailto:}{ykim@culverhouse.ua.edu}.;
\protect \\
H.D. Kwon: Gies College of Business, University of Illinois at Urbana-Champaign;
\protect\href{mailto:}{dhkwon@illinois.edu}.}}
\maketitle
\begin{abstract}
We consider a two-player game of war of attrition under complete information.
It is well- known that this class of games admits equilibria in pure,
as well as mixed strategies, and much of the literature has focused
on the latter. We show that if the players\textquoteright{} payoffs
whilst in \textquotedblleft war\textquotedblright{} vary stochastically
and their exit payoffs are heterogeneous, then the game admits Markov
Perfect equilibria in pure strategies only. This is true irrespective
of the degree of randomness and heterogeneity, thus highlighting the
fragility of mixed-strategy equilibria to a natural perturbation of
the canonical model. In contrast, when the players\textquoteright{}
flow payoffs are deterministic or their exit payoffs are homogeneous,
the game admits equilibria in pure and mixed strategies.
\end{abstract}

\section{Introduction\label{sec:Intro}}

In the classic war of attrition, the first player to quit concedes
a prize to his opponent. Thus, each player trades off the cost associated
with fighting against the value of the prize. These features are common
in many managerial and economic problems. Oligopolists in a declining
industry may bear losses in anticipation of profitability following
a competitor's exit \citep{Ghemawat1985}. For example, the rise of
Amazon in the mid-1990s made the business model of Barnes \& Noble
and Borders obsolete, turning traditional bookselling into a declining
market. As the demand shrank sharply, these two major players at the
time had to cut down slack in their capacities, but each would prefer
its competitors to carry the painful burden of closing stores or exiting
the market altogether \citep{newman2011}. Similarly, the presently
low price of crude oil is often attributed to a war of attrition among
the OPEC allies and non-OPEC rivals such as Russia and the many shale-oil
producers in the United States \citep{Reed2016}. Other examples of
wars of attrition include the provision of public goods \citep{BlissNalebuff1984},
lobbying \citep{Becker1983}, labor disputes \citep{Greenhouse1999},
court of law battles \citep{McAfee2009}, races to dominate a market
\citep{Ghemawat1997}, technology standard races \citep{Bulow1999},
price cycles in oligopolistic collusion \citep{Maskin1988}, all-pay
auctions \citep{KM1997}, and bargaining games \citep{AG2000}.

A central feature of wars of attrition is the waste of valuable resources
(\emph{a.k.a attrition}): there exist times when players would collectively
be better off if one of them quits, but each player strategically
resists quitting in anticipation that his opponent will be the first
to do so. Canonical, complete-information games of war of attrition
typically admit equilibria in both pure and mixed strategies; see,
for example, \citet{Tirole1988}, \citet{FT1996}, and \citet{Levin2004}.
Attrition, however, is featured only in the latter, while the former
are (Pareto) efficient. We study such a two-player model, and show
that if the players' flow payoffs whilst fighting for the prize follow
an irreducible stochastic process and their exit payoffs are heterogeneous,
then the game admits only pure-strategy Markov Perfect equilibria
(hereafter MPE), and this is true irrespective of the degree of uncertainty
or heterogeneity. In other words, our main result shows that an arguably
natural perturbation of the model eradicates all MPE that exhibit
attrition. This result has implications for the modeling choices in
such games, as well as a growing literature that aims to empirically
study strategies in real-world wars of attrition; see, for example,
\citet{Wang2009} and \citet{Takahashi2015}.

In our continuous-time model, two competing oligopolists contemplate
exiting a market. While both firms remain in the market, each receives
the flow payoff that depends on the \emph{market conditions} (\emph{e.g.},
the price of a relevant commodity), which fluctuate according to a
stochastic diffusion process, hereafter the \emph{state}. At every
moment, each firm can exit the market and collect its outside option.
Its rival then obtains a (higher) \emph{winner's} payoff, which depends
on the state at the time of the opponent's exit; \emph{e.g.}, the
net present value of monopoly profits. All payoff-relevant parameters
are common knowledge. The firms may have heterogeneous outside options
but they are otherwise identical. Given that the state follows a Markov
process and mixed strategy equilibria characterized in the literature
are typically stationary (e.g., \citet{Tirole1988}), we focus on
Markov strategies in the main body of the paper, wherein at every
moment, each firm conditions its probability of exit on the current
state.

We first characterize the best response of a firm that anticipates
its rival will never exit, which turns out to be instrumental for
the equilibrium analysis. We show that a firm optimally exits at the
first moment that the state drifts below a threshold. This \emph{single-player-optimal}
threshold is strictly increasing in the firm's outside option; the
better is a firm's outside option, the less it is willing to endure
poor market conditions before exiting.

We present our main result in Section \ref{sec:MPE}. To set the stage,
Proposition \ref{prop:PureMPE} shows that there exists a pure-strategy
MPE in which the firm with the larger outside option exits at the
first moment that the state drifts below its single-player-optimal
threshold. Moreover, if the heterogeneity in outside options is not
too large, then there exists another pure-strategy MPE in which the
firm with the lower outside option exits at the first moment that
the state drifts below its own single-player-optimal threshold. Towards
our main result, we show that in any mixed-strategy MPE, (i) the firms
must randomize between remaining in the market and exiting on a common
set of states, and (ii) each firm exits with nonzero probability if
(and only if) the state is below its single-player-optimal threshold.
However, (i) and (ii) are incompatible with each other if the firms
have heterogeneous outside options, because their single-player-optimal
thresholds differ in that case. Therefore, it follows that no mixed-strategy
MPE exists in that case. We also extend our main result to non-Markovian
Subgame Perfect equilibria subject to a restriction on the firms'
strategies (see Online Appendix \ref{sec:Online-nonMarkov} for details).

\medskip{}

First, this paper contributes to the literature on wars of attrition,
which has received widespread attention since the seminal work of
\citet{Smith1974}. Our model is closest to \citet{Hendricks1988}
and \citet{Murto2004}. The former characterizes equilibria in both
pure and mixed strategies in a war of attrition under complete information
with asymmetric players whose payoffs are deterministic. The latter
considers stochastic payoffs, but restricts attention to pure-strategy
MPE. In contrast, we allow payoffs to vary stochastically, and we
show that if players are heterogeneous, then (subject to a set of
restrictions on strategies) the game admits MPE in pure strategies
\emph{only}.

We also contribute to a literature that contemplates equilibrium selection
in games of war of attrition. This literature has two broad themes.
The first considers games which are \emph{backward-inductible}. For
example, \citet{Ghemawat1985} studies a game with asymmetric players
in which there is a state (that is reached with probability one) at
which both firms have a dominant strategy to exit, while \citet{BS1996}
considers a finite-horizon war of attrition game. In both cases, the
game is shown to have a unique equilibrium in pure strategies. In
the second theme, with a \emph{small} probability, each player \emph{never}
exits. In \citet{FT1986}, players are uncertain about their rivals'
costs of remaining in the market, whereas in \citet{KornhauserEtAl1989},
\citet{Kambe1999}, and \citet{AG2000}, with a small probability,
each player is irrational and never exits. It is shown that the respective
games admit a unique equilibrium. \textcolor{black}{\citet{Myatt2005}
shows that this uniqueness is insensitive to perturbations having
a similar economic interpretation as exit failure.} We complement
this literature by considering a complete-information framework with
rational players, and showing that an arguably natural perturbation
of the canonical model eliminates all mixed-strategy MPE.

\citet{Touzi2002} introduces the concept of mixed strategies in continuous-time
Dynkin games, and proves that the game admits minimax solutions in
mixed strategies. With this notion of mixed strategies, \citet{Strack2016}
investigates a war of attrition with privately observed Brownian
motions, and \citet{Steg2015} characterizes equilibria in both pure
and mixed strategies in a family of continuous-time stochastic timing
games. Whereas these articles consider games with identical players,
we focus on ones with heterogeneous players and show that the set
of equilibria differs drastically. \citet{RS2017} examines mixed-strategy
equilibria in continuous-time stopping games with heterogeneous players,
but focuses on games with pre-emption incentives, whereas ours is
purely one of war of attrition.

\section{Model\label{sec:Model} }

We consider a war of attrition with complete information between two
oligopolists. Time is continuous, and firms discount time at rate
$r>0$. At every moment, each firm decides whether to exit the market.

While both firms remain in the market, each earns a flow profit $\pi(X_{t})$,
where $\pi:\,\mathbb{R}\rightarrow\mathbb{R}$ is continuous and strictly
increasing, and $X_{t}$ is a scalar that captures the \emph{market
conditions} that the firms operate in (\emph{e.g.}, the size of the
market or the price of raw materials). The market conditions fluctuate
according to
\begin{equation}
dX_{t}=\mu(X_{t})dt+\sigma(X_{t})dB_{t}\,,\label{eq:State}
\end{equation}
where $X_{t}$ is defined on $\mathcal{X}:=(\alpha,\beta)\subseteq\mathbb{R}$,
$X_{0}\in\mathcal{X}$, the functions $\mu:\mathcal{X}\rightarrow\mathbb{R}$
and $\sigma:\mathcal{X}\rightarrow\mathbb{R}_{+}$ are Lipschitz continuous,
and $B_{t}$ is a Wiener process.\footnote{Special cases in which $\sigma(\cdot)=0$ have been analyzed extensively
(\citealp{Ghemawat1985}, \citealp{Hendricks1988}, and others). Therefore,
we restrict attention to $\sigma(\cdot)>0$ in the main body of this
paper, and for completeness, we revisit the case in which $\sigma(\cdot)=0$
in Appendix \ref{sec:Symmetric-Deterministic}.}$^{,}$\footnote{The boundary points $\alpha$ and $\beta$ are assumed to be \emph{natural}
\citep[p.18-20]{Borodin1996}; \emph{i.e.}, neither $\alpha$, nor
$\beta$ can be reached by $X_{t}$ in finite time. For example, if
$X_{t}$ is a standard diffusion process, then $\mathcal{X}=\mathbb{R}$.
If $X_{t}$ is a geometric Brownian process, then $\mathcal{X}=(0,\infty)$.} Then because we assume $\sigma(\cdot)>0$ on $\mathcal{X}$, the
functions $\mu(\cdot)$ and $\sigma(\cdot)$ satisfy the local integrability
condition \citep{Arkin2015}, which implies that the diffusion process
$X$ is \emph{regular} in $\mathcal{X}=(\alpha,\beta)$ \citep{KS1991}:
For any $x,y\in\mathcal{X}$, the process $X$ reaches from $x$ to
$y$ in finite time with positive probability; \emph{i.e.}, $X$ is
\emph{irreducible} \citep[p.13]{Borodin1996}. Let $(\Omega,\mathcal{F},\mathbb{P},\{\mathcal{F}_{t}\}_{t\ge0})$
denote the probability space with sample space $\Omega$, $\sigma$-algebra
$\mathcal{F}$, probability measure $\mathbb{P}$, and filtration
$\{\mathcal{F}_{t}\}_{t\ge0}$ that satisfies the \emph{usual} conditions
(p. 172, \citealt{Rogers2000}). We assume that the process $\{B_{t}\}_{t\ge0}$
(or equivalently, $\{X_{t}\}_{t\ge0}$) is progressively measurable
with respect to $\{\mathcal{F}_{t}\}_{t\ge0}$. Throughout the paper,
we let $\mathbb{E}[\cdot]$ denote the expected value with respect
to $\mathbb{P}$.

If firm $i$ chooses to exit at $t$, then it receives its outside
option $l_{i}$, and its opponent, denoted by $-i$, receives $w(X_{t})\in\mathbb{R}$,
the expected payoff associated with being the sole remaining firm;\emph{
e.g.}, the net present value of monopoly profits. In this case, we
say that firm $i$ is the \emph{loser} and firm $-i$ is the \emph{winner}.
We adopt the convention that $l_{1}\leq l_{2}$; \emph{i.e.}, firm
$2$ has a larger outside option than firm $1$. We assume that $w(x)>l_{2}$
for all $x$ so that the winner's reward is always larger than the
loser's. The game ends as soon as a firm exits the market. If both
firms exit at the same moment, then each obtains the outside option
$l_{i}$ or $w(X_{t})$ with probability $1/2$.\footnote{For simplicity, we assume that the firms can differ only in the loser's
exit payoff, $l_{1}$ and $l_{2}$, which is independent of $X$.
In Online Appendix \ref{sec:Structural-Stability}, we show that under
certain conditions, our main result continues to hold if the firms
have heterogeneous discount rates, flow profits, and winner payoffs,
and $l_{i}$ is a function of $X$.}

Finally, we make the following assumptions on the functions $\pi(\cdot)$
and $w(\cdot)$: First, we assume that $\pi(\cdot)$ satisfies the
absolute integrability condition $\mathbb{E}\left[\int_{0}^{\infty}\left|e^{-rt}\pi(X_{t})\right|dt\right]<\infty$,
which ensures that each firm's expected discounted payoff is well-defined
(see \citealp{Alvarez2001}). Second, we assume $w(\cdot)\in C^{2}(\mathcal{X})$
and $w(x)>\mathbb{E}\left[\int_{0}^{t}e^{-rs}\pi(X_{s})ds+\right.$
$\left.e^{-rt}w(X_{t})\,|X_{0}=x\right]$ for all $x\in\mathcal{X}$
and $t$, so that each firm prefers to become the winner sooner rather
than later.\footnote{This assumption is satisfied if and only if $\sigma^{2}(x)w^{\prime\prime}(x)/2+\mu(x)w^{\prime}(x)+\pi(x)>rw(x)$
for all $x\in\mathcal{X}$.} Lastly, we assume that for each $i$, there exists some $x_{ci}\in\mathcal{X}$
such that $\pi(x_{ci})=rl_{i}$, which guarantees the existence of
an optimal exit threshold in the interior of $\mathcal{X}$ (see Lemma
\ref{lemm:Aux_1} and the proof of Lemma \ref{lemm:Opt-Stop} for
details).

\subsection{Markov Strategies \label{subsec:Strategies}}

We assume that both firms employ Markov strategies: At every moment
$t$, each firm chooses (probabilistically) whether to exit based
on the current state $X_{t}$, conditional on the game not having
ended. Formally, each firm $i$ chooses
\begin{enumerate}
\item [i.]a closed subset $E_{i}$ of the state space $\mathcal{X}$ (or
an \emph{exit region}) such that it exits with probability $p_{i}(X_{t})=1$
if $X_{t}\in E_{i}$,
\item [ii.]a non-negative function $\lambda_{i}:\mathcal{X}\rightarrow\mathbb{R}_{+}$
(or an \emph{exit rate}) such that $\lambda_{i}(x)$ represents the
firm's hazard rate of exit when $X_{t}=x$.
\end{enumerate}
Note that we stipulate that $E_{i}$ is a closed set and the exit
probability $p_{i}(x)$ is always $1$. These assumptions can be safely
made because $X$ is a regular diffusion process: For any $x\in\mathcal{X}$,
the hitting times $\tau_{x}^{+}=\inf\{t>0:X_{t}^{x}>x\}$ and $\tau_{x}^{-}=\inf\{t>0:X_{t}^{x}<x\}$
are both 0 almost surely (e.g., \citet[Exercise 3.22, p.312]{RY1991}),
which implies $\tau_{x}=\inf\{t>0:X_{t}^{x}=x\}=0$ almost surely,
i.e., $X$ comes back to the original value indefinitely many times
within any finite time span. This has two implications for Markov
strategies. First, exiting when $X_{t}\in E_{i}$ is indistinguishable
from exiting when $X_{t}\in\text{cl}(E_{i})$. Second, any Markov
strategy in which a firm exits with probability $p_{i}(x)\in(0,1)$
whenever $X_{t}=x$ is indistinguishable from one in which the firm
exits with probability $1$ at the moment $X_{t}=x$.\footnote{Note that the rival's strategy is not relevant to this property of
our Markov strategy. First, in any MPE, only one firm may exit with
positive probability at any $x\in\mathcal{X}$. This is because the
payoff from simultaneous exit $(l_{i}+w(x))/2$ is strictly less than
the winner's payoff $w(x)$. Second, suppose that one firm exits with
positive probability at the hitting time of $x$ while its rival exits
with positive probability at the hitting time of $y$ where $|y-x|=\delta>0$.
Then because of the mentioned property of a regular diffusion process
$X$, if $X$ starts from $x$, it will return to $x$ indefinitely
many times within an arbitrarily small time interval without hitting
$y$, no matter how small $\delta>0$ is.} In the proof of Lemma \ref{lemma:MixedMPE2}, we prove that $E_{i}$
always has an equivalent closed set representation even if it is not
closed.

Throughout this paper, we impose a regularity condition on the function
$\lambda_{i}$. We stipulate that $\lambda_{i}(X_{t})$ is Riemann
integrable over any given time interval $[u,v]$ for $0\leq u<v<\infty$.
Note that firm $i$'s probability of exit within a time interval $[u,v]$
is given by $1-\exp[-\int_{u}^{v}\lambda_{i}(X_{t})dt]$; intuitively,
the Riemann integrability thus ensures that, if the time interval
is discretized, the process of coarse graining of the time interval
does not alter the probability of exit. Conversely, if $\lambda_{i}(X_{t})$
were not Riemann integrable, the lower Riemann summation over time
does not coincide with the upper Riemann summation, which implies
that the continuation time limit is not well-defined. Therefore, it
is natural to impose a condition that $\lambda_{i}(X_{t})$ is Riemann
integrable over any time interval.

We represent firm $i$'s strategy as the pair $a_{i}=(E_{i},\lambda_{i})$,
and $\left\{ a_{1},a_{2}\right\} $ is a \emph{strategy profile}.\footnote{Our definition of a strategy implies that it can be alternatively
expressed as a sum of the absolutely continuous function in time and
the discontinuous jumps without the singularly continuous component.
In Online Appendix \ref{sec:Singular}, we explain how this formulation
of a strategy can be justified.} Intuitively, during any \emph{small} interval $[t,t+dt)$, firm $i$
exits with probability
\[
\rho_{i}(X_{t})=\begin{cases}
1 & \text{ if }X_{t}\in E_{i}\text{ ,}\\
\lambda_{i}(X_{t})dt & \text{ otherwise.}
\end{cases}
\]

If firm $i$ does not exit with probability $1$ at all, we write
$E_{i}=\emptyset$. If it does not exit with a positive hazard rate
(\emph{i.e.,} $\lambda_{i}(X_{t})=0$ for all $t$ almost surely,
hereafter a.s), we write $\lambda_{i}\equiv\mathbf{0}$. Finally,
we say that firm $i$'s strategy is\emph{ pure} if $\lambda_{i}\equiv\mathbf{0}$,
and it is \emph{mixed} otherwise. 

\subsection{Payoffs \label{subsec:Model-Payoff}}

Fix an initial value $X_{0}=x\in\mathcal{X}$ and a strategy profile
$\left\{ a_{i},a_{-i}\right\} $, and define $\tau_{i}:=\inf\{s\geq0:X_{s}\in E_{i}\}$.
Then firm $i$'s exit probability, that is, the probability that firm
$i$ has exited by time $t$, given that $X_{0}=x$, can be written
as
\[
G_{i}(t):=1-(1-\mathbf{1}_{\{t\geq\tau_{i}\}}(t))e^{-\int_{0}^{t}\lambda_{i}(X_{s})ds}\:.
\]
We can define firm $i$'s expected payoff under $\left\{ a_{i},a_{-i}\right\} $
when $X_{0}=x$ as follows:
\begin{align*}
V_{i}(x;a_{i},a_{-i})= & \mathbb{E}\Biggl[\int_{0}^{\infty}\int_{0}^{\infty}\biggl[\int_{0}^{t\land s}e^{-ru}\pi(X_{u})du\\
 & +e^{-r(t\land s)}\Bigl(l_{i}\mathbf{1}_{\{t<s\}}(t)+w(X_{s})\mathbf{1}_{\{t>s\}}(t)+m_{i}(X_{t})\mathbf{1}_{\{t=s\}}(t)\Bigr)\biggr]dG_{-i}(s)dG_{i}(t)|X_{0}=x\Biggr]\,,
\end{align*}
where $m_{i}(x):=(l_{i}+w(x))/2$. The first line represents the
firm's discounted flow payoff until either firm exits. The second
line captures the lump-sum payoff from becoming the winner, becoming
the loser, and exiting simultaneously, respectively.

A strategy profile $\left\{ a_{1}^{*},a_{2}^{*}\right\} $ is a Markov
Perfect equilibrium (MPE) if 
\[
V_{i}(x;a_{i}^{*},a_{-i}^{*})\ge V_{i}(x;a_{i},a_{-i}^{*})
\]
for each firm $i$, every initial value $x$, and every strategy $a_{i}$.

\section{Equilibrium Analysis \label{sec:MPE}}

In Section \ref{subsec:BR}, we characterize the best response of
a firm that anticipates its rival will never exit, which is instrumental
for the equilibrium analysis. We then characterize pure-strategy MPE
in Section \ref{subsec:Pure}, and in Section \ref{subsec:Mixed},
we consider mixed-strategy MPE. In particular, we establish necessary
conditions that any mixed-strategy MPE must satisfy, and our first
main result follows immediately: The game has no mixed-strategy MPE
if the firms have heterogeneous exit payoffs (\emph{i.e.,} $l_{1}\neq l_{2}$).

\subsection{Best Response to a Firm which Never Exits\label{subsec:BR}}

We characterize firm $i$'s best response assuming that its rival
never exits;\emph{ i.e.,} the best response to $a_{-i}=\left(\emptyset,\mathbf{0}\right)$.
In this case, the firm's best response can be determined by solving
a single-player optimal stopping problem as in \citet{Alvarez2001}.
Because $X$ is a time-homogeneous process and the time horizon is
infinite, it is without loss of generality to restrict attention to
pure strategies such that $\lambda_{i}=\mathbf{0}$, and the firm's
expected payoff at $t$ depends solely on the current value of the
state $x=X_{t}$.\footnote{Note that any strategy with $\lambda_{i}\neq\mathbf{0}$ \emph{mixes}
pure strategies. Thus, if a strategy with $\lambda_{i}\neq\mathbf{0}$
is a best response to $a_{-i}=\left(\emptyset,\mathbf{0}\right)$,
then there must exist more than one stopping times that are solutions
of the single-player optimal stopping problem, \eqref{eq:Opt-Stop}.
However, this optimal stopping problem admits a unique solution, which
is a hitting time, given in Lemma \ref{lemm:Opt-Stop} (\emph{e.g.},
see \citet{Alvarez2001,Arkin2015}).} Thus, it can be expressed as 
\begin{equation}
\sup_{\tau_{i}\geq t}\,\mathbb{E}\biggl[\int_{t}^{\tau_{i}}e^{-r(s-t)}\pi(X_{s})ds+e^{-r(\tau_{i}-t)}l_{i}|X_{t}=x\biggl]\,.\label{eq:Opt-Stop}
\end{equation}
Using Proposition 2 in \citet{Alvarez2001}, we can characterize the
firm's optimal exit region as follows.\begin{lemma} \label{lemm:Opt-Stop}
Suppose firm $-i$ never exits. There exists a unique threshold $\theta_{i}^{*}$
such that $E_{i}^{*}=\left\{ X_{t}\leq\theta_{i}^{*}\right\} $ is
optimal for firm $i$, that is, firm $i$ optimally exits whenever
$X_{t}\leq\theta_{i}^{*}$. If $l_{1}<l_{2}$, then $\theta_{1}^{*}<\theta_{2}^{*}$.\end{lemma}

The proof is relegated to Online Appendix \ref{sec:Lemma1}. Intuitively,
a firm's value of remaining in the market decreases as the market
conditions deteriorate, and once they become sufficiently poor, the
firm is better off exiting and collecting its outside option. As the
firms earn identical flow payoffs while in the market, the firm with
the higher outside option optimally exits at a higher threshold.

\subsection{Pure-strategy MPE \label{subsec:Pure}}

The following proposition shows that there is a pure-strategy MPE
in which firm $2$ exits at the first moment that $X_{t}$ drifts
below $\theta_{2}^{*}$, and firm $1$ never exits. Moreover, if the
firms are not too heterogeneous, there is another pure-strategy MPE
in which firm $1$ exits at the first moment that $X_{t}\leq\theta_{1}^{*}$
and firm $2$ never exits. \begin{proposition}\label{prop:PureMPE}
~

\noindent (i) The strategy profile $\left\{ a_{1},a_{2}\right\} =\left\{ \left(\emptyset,\mathbf{0}\right),\left(E_{2}^{*},\mathbf{0}\right)\right\} $
is a pure-strategy MPE, where $E_{i}^{*}=(\alpha,\theta_{i}^{*}]$
and $\theta_{i}^{*}$ is given in Lemma \ref{lemm:Opt-Stop}. 

\noindent (ii) There exists a threshold $\kappa>0$ that is independent
of $l_{1}$ such that $\left\{ a_{1},a_{2}\right\} =\left\{ (E_{1}^{*},\mathbf{0}),(\emptyset,\mathbf{0})\right\} $
is also a pure-strategy MPE if $|l_{2}-l_{1}|<\kappa$. 

\noindent (iii) If $X_{0}\geq\max\{\theta_{1}^{*},\theta_{2}^{*}\}$,
every pure-strategy MPE is payoff-/outcome-equivalent to one of the
above.\end{proposition} 

The proof is provided in Appendix \ref{sec:Proofs}. If firm $i$
expects its rival to never exit, then by Lemma \ref{lemm:Opt-Stop},
it will optimally exit at the first time such that $X_{t}\leq\theta_{i}^{*}$.
Therefore, it suffices to show that if firm $i$ employs the exit
region $E_{i}^{*}$, then its opponent's best response is to never
exit.

Suppose that firm $1$ expects its rival to exit at the first moment
that $X_{t}\leq\theta_{2}^{*}$. Recall that firm $2$ has a better
outside option than firm $1$ (\emph{i.e.}, $l_{2}\geq l_{1}$), so
by Lemma \ref{lemm:Opt-Stop}, $\theta_{1}^{*}\leq\theta_{2}^{*}$,
which implies that firm $1$ has no incentive to exit until at least
$X_{t}\leq\theta_{1}^{*}$. Therefore, firm $1$ expects that the
game will end before the state hits $\theta_{1}^{*}$, and hence the
strategy of never exiting is incentive compatible. If instead firm
$2$ anticipates that its rival chooses $(E_{1}^{*},\mathbf{0})$,
then the strategy $a_{2}=(\emptyset,\mathbf{0})$ is incentive compatible
as long as it does not need to wait too long until $X_{t}$ hits $\theta_{1}^{*}$
and firm $1$ exits. As a result, never exiting is a best response
for firm $2$ as long as $\left|l_{2}-l_{1}\right|$, and hence $\theta_{2}^{*}-\theta_{1}^{*}$
is not too large.

If $X_{0}<\max\{\theta_{1}^{*},\theta_{2}^{*}\}$ and $\sigma(\cdot)$
is sufficiently large, as shown in Proposition 5 in \citet{Murto2004},
there may also exist pure-strategy MPE with multiple exit thresholds.
As such equilibria do not affect our analysis of mixed-strategy equilibria,
we do not consider them here. Finally, because along the path of any
equilibrium characterized in Proposition \ref{prop:PureMPE}, at most
one player resists exiting below his single-player optimal threshold,
it follows that both equilibria are Pareto-efficient. 

\subsection{Mixed-strategy MPE \label{subsec:Mixed}}

We now consider mixed-strategy MPE. First, we define the support of
firm $i$'s mixed strategy as the subset of the state space in which
firm $i$ randomizes between remaining in the market and exiting,
\begin{equation}
\Gamma_{i}=\{x\in\mathcal{X}:\lambda_{i}(x)>0\}\:.\label{eq:supp-lambda}
\end{equation}
The support $\Gamma_{i}$ can be represented as a union of open intervals
in $\mathcal{X}$ because of the regularity condition imposed in Section
\ref{subsec:Strategies} that $\lambda_{i}(X_{t})$ is Riemann integrable
over time intervals (see Lemma \ref{lemm:Aux_2} for the proof of
this statement). Intuitively, this property of $\Gamma_{i}$ implies
that whenever $X_{t}\in\Gamma_{i}$, firm $i$ exits with probability
$\lambda_{i}(X_{t})dt>0$ during the time interval $[t,t+dt)$ and
continues to randomize its decision until $\tau_{\Gamma_{i}}:=\inf\{s\ge t:X_{s}\not\in\Gamma_{i}\}$.

The following lemma shows that the firms' mixed strategies must
have common support and neither firm exits with probability $1$.

\begin{lemma} \label{lemma:MixedMPE2} Suppose that $\sigma(\cdot)>0$,
and $\left\{ a_{1},a_{2}\right\} $ constitutes a mixed-strategy MPE.
Then the firms' mixed strategies have common support $\Gamma=(\alpha,\theta_{1}^{*})=(\alpha,\theta_{2}^{*})$,
where $\theta_{i}^{*}$ is given in Lemma \ref{lemm:Opt-Stop}, and
$E_{1}=E_{2}=\emptyset$.\end{lemma}

We give a sketch of the proof below, while the formal proof is relegated
to Appendix \ref{sec:Proofs}. At any time $t$ such that $X_{t}\in\Gamma_{i}$,
firm $i$ must be indifferent between exiting immediately and remaining
in the market, which implies that its expected payoff must be equal
to its outside option, that is,
\begin{equation}
l_{i}=\lambda_{-i}(X_{t})dt\,w(X_{t})+\left(1-\lambda_{-i}(X_{t})dt\right)\left[\pi(X_{t})dt+(1-rdt)l_{i}\right]\,.\label{eq:Indifference}
\end{equation}
The left-hand-side of (\ref{eq:Indifference}) represents firm $i$'s
payoff if it exits at $t$, while the right-hand side represents its
payoff if it remains. To be specific, with probability $\lambda_{-i}(X_{t})dt$,
it receives the winner's payoff, $w(X_{t})$, whereas with the complementary
probability, it earns the flow payoff $\pi(X_{t})$ during $(t,t+dt)$,
and its (discounted) continuation profit, $l_{i}$, at $t+dt$.\footnote{We ignore the event that both firms exit simultaneously. As the proof
shows, this is an innocuous simplification.} Thus, firm $-i$'s exit rate must satisfy
\begin{equation}
\lambda_{-i}(X_{t})=\frac{rl_{i}-\pi(X_{t})}{w(X_{t})-l_{i}}\,.\label{eq:ExitRate}
\end{equation}
Note that $\pi(x)<rl_{i}$ for any $x\in\Gamma_{i}$.\footnote{If $\pi(X_{t})>rl_{i}$, then the right-hand-side of (\ref{eq:Indifference})
is strictly larger than $l_{i}$, so firm $i$ strictly prefers to
remain in the market regardless of its rival's strategy.} We shall now argue that $\Gamma_{1}=\Gamma_{2}$. Towards a contradiction,
suppose that there exists a non-empty interval that is a subset of
$\Gamma_{i}$ but not of $\Gamma_{-i}$. Then for any $x$ in that
interval, $\pi(x)<rl_{i}$ and $\lambda_{-i}(x)=0$, because by assumption,
$x\in\Gamma_{i}$ and $x\notin\Gamma_{-i}$, respectively. This implies
that the right-hand-side of (\ref{eq:Indifference}) is strictly smaller
than $l_{i}$, so firm $i$ strictly prefers to exit, which contradicts
that $x\in\Gamma_{i}$. Hence, we conclude that $\Gamma_{i}\backslash\Gamma_{-i}$
is empty, and so $\Gamma_{1}=\Gamma_{2}$.

Next, recall that even if firm $i$ anticipates that its rival will
never exit, it is unwilling to exit until $X_{t}$ hits $\theta_{i}^{*}$.
Hence, if this firm expects its rival to exit with positive probability,
then \emph{ceteris paribus}, this \emph{decreases} its incentive to
exit. Consequently, firm $i$ \emph{always} strictly prefers to remain
in the market whenever $X_{t}>\theta_{i}^{*}$, which, together with
the fact that $\Gamma$ is open, implies that $\Gamma\subseteq(\alpha,\theta_{i}^{*})$.

We now argue that in a mixed-strategy MPE, neither firm exits with
probability $1$, that is, $E_{1}=E_{2}=\emptyset$. Towards a contradiction,
suppose that $E_{i}\neq\emptyset$ for some $i$. Because exiting
at any $X_{t}>\theta_{i}^{*}$ is a strictly dominated strategy for
firm $i$, it must be the case that $E_{i}\subseteq(\alpha,\theta_{i}^{*}]$.
Moreover, because firm $-i$ strictly prefers to remain in the market
when $X_{t}$ is sufficiently close to $E_{i}$ (anticipating that
firm $i$ will soon exit with probability one), it must also be the
case that $E_{i}$ and $\Gamma$ are disjoint and separated by a non-empty
interval $(c,d)$ wherein neither firm exits. As both $E_{i}$ and
$\Gamma$ are subsets of $(\alpha,\theta_{i}^{*}]$, so must be $(c,d)$.
Then because $\pi(x)<rl_{i}$ for any $x\in(c,d)\subseteq(\alpha,\theta_{i}^{*}]$,
firm $i$ strictly prefers to exit instantaneously if $X_{t}\in(c,d)$,
instead of waiting until the state hits $E_{i}$ or $\Gamma$, contradicting
the premise that $a_{i}$ is a best response to $a_{-i}$. Hence,
we conclude that $E_{i}=\emptyset$.

We have already argued that $\Gamma\subseteq(\alpha,\theta_{i}^{*})$.
It remains to argue that this inclusion is an equality. Suppose that
$\Gamma=(\alpha,\theta)$ for some $\theta<\theta_{i}^{*}$. Because
firm $-i$ does not exit at any $X_{t}>\theta$, firm $i$'s expected
payoff at any $x\in(\theta,\theta_{i}^{*})$ from exiting at the first
time that $X_{t}\leq\theta$ is strictly less than $l_{i}$ by Lemma
\ref{lemm:Opt-Stop}, which implies that this firm strictly prefers
to exit instantaneously--a contradiction.\medskip{}

Because $\theta_{1}^{*}<\theta_{2}^{*}$ whenever $l_{1}<l_{2}$ by
Lemma \ref{lemm:Opt-Stop}, we have the following immediate implication.

\begin{theorem} \label{thm:Stoch-MixedMPE} Suppose that $\sigma(\cdot)>0$
and $l_{1}<l_{2}$. Then the game admits no mixed-strategy MPE. \end{theorem}

This theorem, together with Proposition \ref{prop:PureMPE}, shows
that if there is even a small amount of uncertainty about the payoff
from remaining in the market and the firms are even slightly heterogeneous,
then none of the MPE feature \emph{attrition}, \emph{i.e., }both firms
resisting exit below their single-player-optimal thresholds. 

Both conditions in Theorem \ref{thm:Stoch-MixedMPE} are necessary
to eliminate mixed-strategy MPE. If the firms are homogeneous (\emph{i.e.,}
$l_{1}=l_{2}$) or payoffs are deterministic (\emph{i.e.,} $\sigma(\cdot)\equiv0$),
then as shown in \citet{Steg2015} and \citet{Hendricks1988}, respectively,
and, for completeness, as we show in Appendix \ref{sec:Symmetric-Deterministic},
the game admits MPE in both pure and mixed strategies. 

The key driver behind this result is that the state follows an irreducible
stochastic process. If the state is deterministic, as in \citet{Hendricks1988}
for example, then in any mixed-strategy equilibrium, the firm with
the smaller outside option (firm 1) exits with a positive probability
when the state hits its single-player optimal threshold, and from
then onward, during every interval of length $dt$, each firm exits
with some probability that is proportional to $dt$ and makes the
other firm indifferent between exiting and not. If the state is stochastic,
any hypothetical candidate mixed-strategy MPE must also be of the
form described above. However, when the state is just below that threshold,
due to the irreducibility of the state, with some likelihood, it will
hit that threshold in short order and firm 1 will exit. As a result,
firm 2 strictly prefers to not exit near that threshold, which in
turn leads firm 1 to strictly prefer to exit, leading to a pure-strategy
MPE. \medskip{}

A natural concern is the restrictiveness of Markov strategies. Toward
investigating this, we consider the possibility that firms condition
their exit decision on the history. Under a set of restrictions on
the firms' strategies, we show that with heterogeneous outside options,
the game admits no mixed-strategy Subgame Perfect equilibria in Online
Appendix \ref{sec:Online-nonMarkov}. To elaborate on the restrictions,
note that a non-Markovian strategy consists of (i) a set of stopping
times at which the firm exits with positive probability, and (ii)
an \emph{exit rate} function, specifying the probability that the
firm exits during an interval $(t,t+dt)$, which depends on the history
of the state up to $t$. Then our restriction imposes that the firms
exit with positive probability at no more than finitely-many stopping
times, they do not exit with probability one, and the exit rate function
satisfies a regularity condition analogous to that in Markov strategies.
We remark that these restrictions are satisfied by the strategies
in the mixed-strategy SPE that appear in the literature when the state
evolves deterministically \citep{Hendricks1988}, or the firms are
homogeneous \citep{Steg2015}, or both \citep{Tirole1988}.\footnote{If the state $X$ always drifts downward, then there is a one-to-one
correspondence between the state $X$ and time $t$, which implies
that every SPE is also an MPE, i.e., the class of SPEs coincides with
that of MPEs.}

\section{Discussion}

We consider a two-player war of attrition under complete information.
Our main result shows that if the players are heterogeneous and their
flow payoffs whilst in war follow a diffusion process, then the game
admits no mixed-strategy MPE. We also extend this non-existence result
to a class of Subgame Perfect equilibria (subject to a set of restrictions),
where the players' strategies can depend on the entire history. Because
the pure-strategy equilibria are Pareto-efficient, these results indicate
that an arguably natural perturbation of the canonical model eradicates
equilibria that possess a central feature of wars of attrition---the
waste of valuable resources, suggesting that the complete-information
model may be unsuitable for studying this class of problems.

Much of the recent theoretical and empirical literature on wars of
attrition has focused on asymmetric-information models (\emph{e.g.,}
\citet{Myatt2005}, \citet{Wang2009} and \citet{Takahashi2015}),
which admit equilibria that feature attrition. However, as these equilibria
are often obtained by purifying mixed strategies in the complete-information
game, it is an open question whether there exist any equilibria that
feature attrition in an incomplete-information counterpart of our
model. We view this as a promising avenue for future research.\newpage{}

\bibliographystyle{INFORMS2011}
\bibliography{WA}

\appendix

\section{Mixed Strategy MPE in two Special Cases ($l_{1}=l_{2}$ or $\sigma(\cdot)\equiv0$)
\label{sec:Symmetric-Deterministic}}

Recall that if\emph{ $\sigma(\cdot)>0$} and $l_{1}<l_{2}$, then
the game admits no mixed-strategy MPE (Theorem \ref{thm:Stoch-MixedMPE}).
In this section, we show that if either of these conditions is removed,
a mixed-strategy MPE does exist.

First, let us consider the case in which the firms are homogeneous
(\emph{i.e.}, $l_{1}=l_{2}$) and \emph{$\sigma(\cdot)>0$}. It follows
from Lemma \ref{lemm:Opt-Stop} that $\theta_{1}^{*}=\theta_{2}^{*}$.
Following \citet{Steg2015}, it is easy to show that the strategies
$a_{1}=\left(\emptyset,\lambda_{1}(\cdot)\right)$ and $a_{2}=\left(\emptyset,\lambda_{2}(\cdot)\right)$,
where 
\begin{equation}
\lambda_{i}(x):=\mathbb{I}_{\left\{ x\leq\theta_{1}^{*}\right\} }\frac{rl_{-i}-\pi(x)}{w(x)-l_{-i}}\,,\label{eq:lambda_i}
\end{equation}
 constitute a mixed strategy MPE. Note that $\lambda_{1}=\lambda_{2}$
in this case because $l_{1}=l_{2}$.

Next, we consider the case in which $X$ evolves deterministically
(\emph{i.e.}, $\sigma(\cdot)\equiv0$) and $l_{1}<l_{2}$. From Lemma
\ref{lemm:Opt-Stop}, we have that $\theta_{1}^{*}<\theta_{2}^{*}$.
Following \citet{Hendricks1988}, we let $\mu(\cdot)\leq0$, i.e.,
the market condition always goes down over time. Then unlike the case
in which $\sigma(\cdot)>0$ considered in Section \ref{subsec:Mixed},
when $X$ is deterministic, it is not without loss of generality to
assume that $p_{i}(x)=1$, i.e., $p_{i}(x)<1$ is possible in Markov
strategies. Let $E_{1}^{*}=(\alpha,\theta_{1}^{*}]$ with $p_{1}(x)=q_{1}\in(0,1)$
for all $x\in E_{1}^{*}$. Consider the strategies $a_{1}=\left(E_{1}^{*}(p_{1}),\lambda_{1}(\cdot)\right)$
and $a_{2}=\left(\emptyset,\lambda_{2}(\cdot)\right)$, where $E_{1}^{*}(p_{1})$
indicates the exit probability $p_{1}$ whenever $X_{t}\in E_{1}^{*}$
and $\lambda_{i}$ is given in \eqref{eq:lambda_i}. Because $X$
always moves downwards, the strategies defined above are Markov. Using
similar arguments to \citet{Hendricks1988}, one can show that if
$\left|l_{1}-l_{2}\right|$ is not too large, then there exists a
$q_{1}\in(0,1)$ such that $(a_{1},a_{2})$ constitutes a mixed strategy
MPE.

\section{Mathematical Supplement\label{sec:Math-Supplement}}

This section provides supplementary lemmas that are used to prove
the results (lemmas and propositions) in the body of the manuscript.

We first define the following functions that will be used later.
\begin{align}
R(x) & :=\mathbb{E}^{x}\biggl[\int_{0}^{\infty}\pi(X_{t})e^{-rt}dt\biggr]\:,\label{eq:R}\\
\beta_{i}(x) & :=\frac{l_{i}-R(x)}{\phi(x)}\:,\label{eq:Beta_i}
\end{align}
where $\phi:\mathcal{X}\rightarrow\mathbb{R}$ satisfies the differential
equation\footnote{This second-order linear ordinary differential equation (ODE) always
has two linearly independent fundamental solutions, one of which is
monotonically decreasing (see \citealp[p.319]{Alvarez2001}). Note
that if $f(\cdot)$ solves this equation, then so does $cf(\cdot)$
for any constant $c\in\mathbb{R}$ because it is a homogeneous equation.
Hence, we can always find the one which is always positive.} $\frac{1}{2}\sigma^{2}(x)\phi^{''}(x)+\mu(x)\phi^{'}(x)-r\phi(x)=0$
with the properties of $\phi(\cdot)>0$ and $\phi^{'}(\cdot)<0$.
The function $R(\cdot)$ is well-defined because we assume that $\pi(\cdot)$
satisfies the absolute integrability condition in Section \ref{sec:Model}.
The following lemma establishes some properties of the function $\beta_{i}$.
This lemma will be used to prove Lemma \ref{lemm:Opt-Stop} and Proposition
\ref{prop:PureMPE}.

\begin{lemma} \label{lemm:Aux_1} The function $\beta_{i}(x)$ has
a unique interior maximum at $\theta_{i}^{*}\le x_{ci}$ where $\pi(x_{ci})=rl_{i}$.
Furthermore, $\beta_{i}^{'}(x)>0$ for $x<\theta_{i}^{*}$ and $\beta_{i}^{'}(x)<0$
for $x>\theta_{i}^{*}$. \end{lemma}\textbf{Proof} \textbf{of Lemma
\ref{lemm:Aux_1}}: To prove this lemma, it is enough to examine the
behavior of the first derivative of $\beta_{i}(x)=[l_{i}-R(x)]/\phi(x)$.

According to the theory of diffusive processes \citep[p.319]{Alvarez2001},
the function $R(\cdot)$, given in \eqref{eq:R}, can be expressed
as
\begin{align}
R(x) & =\frac{\phi(x)}{B}\int_{a}^{x}\psi(y)\pi(y)m^{'}(y)dy+\frac{\psi(x)}{B}\int_{x}^{b}\phi(y)\pi(y)m^{'}(y)dy\:.\label{eq:Pf-LemmaA1-R}
\end{align}
Here, $a$ and $b$ are the two boundaries of the state space $\mathcal{X}$,
$\psi(\cdot)$ and $\phi(\cdot)$ are the increasing and decreasing
fundamental solutions to the differential equation $\frac{1}{2}\sigma^{2}(x)f^{''}(x)+\mu(x)f^{'}(x)-rf(x)=0$,
$B=[\psi^{'}(x)\phi(x)-\psi(x)\phi^{'}(x)]/S^{'}(x)$ is the \emph{constant}
Wronskian determinant of $\psi(\cdot)$ and $\phi(\cdot)$, $S^{'}(x)=\exp(-\int2\mu(x)/\sigma^{2}(x)dx)$
is the density of the scale function of $X$, and $m^{'}(y)=2/[\sigma^{2}(y)S^{'}(y)]$
is the density of the speed measure of $X$.

By virtue of \eqref{eq:Pf-LemmaA1-R}, differentiation of $R(x)$
with respect to $x$ leads to
\begin{align}
R^{'}(x)\phi(x)-R(x)\phi^{'}(x) & =S^{'}(x)\int_{x}^{b}\phi(y)\pi(y)m^{'}(y)dy\:.\label{eq:Pf-LemmaA1-Rphi}
\end{align}
Moreover, because $l_{i}=\mathbb{E}^{x}[\int_{0}^{\infty}rl_{i}e^{-rt}dt]$,
we can write
\begin{align}
R(x)-l_{i} & =\mathbb{E}^{x}\biggl[\int_{0}^{\infty}[\pi(X_{t})-rl_{i}]e^{-rt}dt\biggr]\:,\label{eq:Pf-LammaA1-RL}
\end{align}
which implies that we can treat the functional $R(x)-l_{i}$ as the
expected cumulative present value of a flow payoff $\pi(\cdot)-rl_{i}$.
Combining \eqref{eq:Pf-LemmaA1-Rphi} and \eqref{eq:Pf-LammaA1-RL},
therefore, we obtain
\begin{align}
\beta_{i}^{'}(x) & =-\frac{R^{'}(x)\phi(x)-[R(x)-l_{i}]\phi^{'}(x)}{\phi^{2}(x)}=-\frac{S^{'}(x)}{\phi^{2}(x)}\int_{x}^{b}\phi(y)[\pi(y)-rl_{i}]m^{'}(y)dy\:.\label{eq:Pf-LemmaA1-Betap}
\end{align}

Now, because $\pi(\cdot)$ is strictly increasing and $\pi(x_{ci})=rl_{i}$,
it must be the case that $\pi(x)<rl_{i}$ for $x<x_{ci}$ and $\pi(x)>rl_{i}$
for $x>x_{ci}$. Thus, $\beta_{i}^{'}(x)<0$ for all $x>x_{ci}$.
Note also that if $x<K<x_{ci}$, then
\begin{align*}
\int_{x}^{b}\phi(y)[\pi(y)-rl_{i}]m^{'}(y)dy & =\int_{x}^{K}\phi(y)[\pi(y)-rl_{i}]m^{'}(y)dy+\int_{K}^{b}\phi(y)[\pi(y)-rl_{i}]m^{'}(y)dy\\
 & \leq\frac{[\pi(K)-rl_{i}]}{r}\biggl(\frac{\phi^{'}(K)}{S^{'}(K)}-\frac{\phi^{'}(x)}{S^{'}(x)}\biggr)+\int_{K}^{b}\phi(y)[\pi(y)-rl_{i}]m^{'}(y)dy\rightarrow-\infty\:,
\end{align*}
as $x\downarrow a$ because $a$ is a natural boundary, which implies
that $\lim_{x\downarrow a}\beta_{i}^{'}(x)=\infty$. Here we use $\phi^{'}(x)<0$
and $\pi(x)<\pi(K)<rl_{i}$ for $x<K$. It thus follows that $\beta_{i}^{'}(\theta_{i}^{*})=0$
for some $\theta_{i}^{*}\leq x_{ci}$, which implies that $\int_{\theta_{i}^{*}}^{b}\phi(y)[\pi(y)-rl_{i}]m^{'}(y)dy=0$
because $S^{'}(x)>0$ and $\phi(x)>0$ in \eqref{eq:Pf-LemmaA1-Betap}.
Moreover, note that $\int_{x}^{b}\phi(y)[\pi(y)-rl_{i}]m^{'}(y)dy$
is increasing in $x<x_{ci}$ because $\pi(y)<rl_{i}$ for $\forall y<x_{ci}$,
thus yielding $\int_{x}^{b}\phi(y)[\pi(y)-rl_{i}]m^{'}(y)dy<0$ if
$x<\theta_{i}^{*}\leq x_{ci}$ and $\int_{x}^{b}\phi(y)[\pi(y)-rl_{i}]m^{'}(y)dy>0$
if $\theta_{i}^{*}<x\leq x_{ci}$. Combining this with \eqref{eq:Pf-LemmaA1-Betap},
we obtain the unique existence of $\theta_{i}^{*}$ such that $\beta_{i}^{'}(x)>0$
for $\forall x<\theta_{i}^{*}$ and $\beta_{i}^{'}(x)<0$ for $\forall x>\theta_{i}^{*}$,
which completes the proof. \eproof

\medskip{}

Recall that we define the support $\Gamma_{i}=\{x\in\mathcal{X}:\lambda_{i}(x)>0\}$
for the hazard rate function $\lambda_{i}$ where we impose the regularity
condition that $\lambda_{i}(X_{t})$ is Riemann integrable over time
intervals. In the following lemma, we prove that the regularity condition
on $\lambda_{i}$ implies $\Gamma_{i}=\{x\in\mathcal{X}:\lambda_{i}(x)>0\}$
is a union of open intervals, hence it is an open set. This lemma
will be used to prove Lemma \ref{lemma:MixedMPE2}.

\begin{lemma} \label{lemm:Aux_2} If $\lambda_{i}(X_{t})$ is Riemann
integrable over any bounded time interval $[u,v]$, then $\Gamma_{i}$
can be represented as a union of open intervals in the state space
$\mathcal{X}$. \end{lemma} \textbf{Proof} \textbf{of Lemma \ref{lemm:Aux_2}:}
For notational simplicity, we let $\lambda_{i,t}:=\lambda_{i}(X_{t})$.
For any given bounded time interval $[u,v]$, if $\lambda_{i,t}$
is Riemann integrable for a given sample path of $X_{u\leq t\leq v}$,
then $\lambda_{i,t}$ must be bounded within the interval $[u,v]$.
Then by virtue of Theorem 2.28 of \citealt{Folland1999}, the set
$D_{i}:=\{t\in[u,v]:\lambda_{i,t}\;\text{is discontinuous}\}$ has
Lebesgue measure zero. Hence, we can construct an equivalent version
of $\lambda_{i,t}$ where $\lambda_{i,t}=0$ whenever $t\in D_{i}$;
because $\int_{D_{i}}\lambda_{i,t}dt=0$ even if $\lambda_{i,t}>0$
for all $t\in D_{i}$, this transformation does not affect either
the firm's payoff or the outcome of the game. By construction, the
transformed process $\{\lambda_{i,t}\}$ is continuous for all $t$
at which $\lambda_{i,t}>0$. Therefore, $C_{i}:=\{t\in[u,v]:\lambda_{i,t}>0\}$
is an open set of the real line $\mathbb{R}$, and hence, a countable
union of disjoint open intervals of time.

Next, because $X_{t}$ is a continuous process in time $t$, the set
$\{X_{t}:t\in C_{i}\}$, which is the mapping of the set $C_{i}$
into the state space via $X$, is a union of disjoint intervals, open
or closed; it may contain point sets, but without loss of generality,
we may ignore them because these possibilities are zero-probability
events. Note that $\{X_{t}:t\in C_{i}\}\subseteq\Gamma_{i}$ by the
definition of $C_{i}$ and $\Gamma_{i}$. Hence, $\Gamma_{i}$ is
a union of disjoint intervals in the state space $\mathcal{X}$ because
otherwise there would be some time interval $[u,v]$ and some sample
path $X$ such that $C_{i}$ cannot be represented as a union of time
intervals.

As a final step, we can equivalently represent $\Gamma_{i}$ as a
union of disjoint \emph{open} intervals by taking its interior (i.e.,
removing all the boundary points of $\Gamma_{i}$, if any) without
affecting the equilibrium payoffs because $X$ is a regular diffusion
process: Let $B_{i}$ be the set of all boundary points of $\Gamma_{i}$.
Then because $\Gamma_{i}$ is a union of disjoint intervals in the
real line $\mathbb{R}$, the cardinality of $B_{i}$ is at most countably
infinite, which implies that $B_{i}$ has Lebesgue measure zero. Hence,
we have $\int_{0}^{t}\mathbf{1}_{\{x\in B_{i}\}}(X_{s})\lambda_{i}(X_{s})ds=0$
because $X$ is a regular diffusion process and neither strategy nor
associated payoff is affected by the removal of any of the boundary
points of $\Gamma_{i}$.\eproof

\section{Proofs\label{sec:Proofs}}

\noindent \textbf{Proof} \textbf{of Proposition \ref{prop:PureMPE}}: 

\noindent (i) Define $\tau_{i}^{*}:=\inf\{t\geq0:X_{t}\in E_{i}^{*}\}$,
$i\in\{1,2\}$, where $E_{i}^{*}=(\alpha,\theta_{i}^{*}]$ is given
in Section \ref{subsec:Pure}. For expositional convenience, we also
let $\tau(E_{i}):=\left(E_{i},\mathbf{0}\right)$ for each $i\in\{1,2\}$.
We first prove that $\{a_{1},a_{2}\}=\{\tau(\emptyset),\tau(E_{2}^{*})\}$
is an MPE. Because it is shown in Lemma \ref{lemm:Opt-Stop} that
$a_{2}=\tau(E_{2}^{*})$ is firm $2$'s best response to $a_{1}=\tau(\emptyset)$,
it only remains to prove that $a_{1}=\tau(\emptyset)$ is also firm
$1$'s best response to $a_{2}=\tau(E_{2}^{*})$.

Let $\tau(E_{1})$ be firm $1$'s best response to $\tau(E_{2}^{*})$
and $V_{W1}^{*}(x):=V_{1}(x;\tau(E_{1}),\tau(E_{2}^{*}))$ be the
corresponding payoff to firm $1$. We let $C_{1}=\mathcal{X}\backslash E_{1}$
denote the continuation region associated with the strategy $\tau(E_{1})$.

First, we show that $E_{1}\cap(\theta_{2}^{*},\infty)=\emptyset$.
Toward a contradiction, suppose this is not the case. Then pick some
$x\in E_{1}\cap(\theta_{2}^{*},\infty)$ and observe that $V_{W1}^{*}(x)=l_{1}$
due to $x\in E_{1}$. However, 
\begin{align*}
V_{W1}^{*}(x)\geq V_{1}(x;\tau(\emptyset),\tau(E_{2}^{*})) & =\mathbb{E}^{x}\biggl[\int_{0}^{\tau_{2}^{*}}\pi(X_{t})e^{-rt}dt+w(X_{\tau_{2}^{*}}^{x})e^{-r\tau_{2}^{*}}\biggr]=R(x)+\biggl[\frac{w(\theta_{2}^{*})-R(\theta_{2}^{*})}{\phi(\theta_{2}^{*})}\biggr]\phi(x)\\
 & >R(x)+\biggl[\frac{l_{1}-R(\theta_{2}^{*})}{\phi(\theta_{2}^{*})}\biggr]\phi(x)=R(x)+\beta_{1}(\theta_{2}^{*})\phi(x)>R(x)+\beta_{1}(x)\phi(x)=l_{1}\:,
\end{align*}
where the first inequality follows because $w(X_{\tau_{2}^{*}}^{x})=w(\theta_{2}^{*})>l_{1}$
and $\mathbb{E}^{x}[e^{-r\tau_{2}^{*}}]=\phi(x)/\phi(\theta_{2}^{*})$
for $x>\theta_{2}^{*}$, and the second inequality holds because $x>\theta_{2}^{*}>\theta_{1}^{*}$
and $\beta_{1}^{'}(x)<0$ for $x>\theta_{1}^{*}$ by Lemma \ref{lemm:Aux_1}.
This establishes the contradiction.

Second, we also prove that $E_{1}\cap(-\infty,\theta_{2}^{*}]=\emptyset$.
Towards a contradiction, suppose this is not the case. Then we can
pick some $x\in E_{1}\cap(-\infty,\theta_{2}^{*}]$ such that $V_{W1}^{*}(x)=m_{1}(x)$
because $\tau_{2}^{*}=\inf\{t\geq0:X_{t}\in E_{2}^{*}\}$. However,
\begin{align*}
V_{W1}^{*}(x)\geq V_{1}(x;\tau(\emptyset),\tau(E_{2}^{*})) & =\mathbb{E}^{x}\biggl[\int_{0}^{\tau_{2}^{*}}\pi(X_{t})e^{-rt}dt+w(X_{\tau_{2}^{*}}^{x})e^{-r\tau_{2}^{*}}\biggr]=w(x)>m_{1}(x)\:,
\end{align*}
where the second equality uses that $\tau_{2}^{*}=0$ when $X_{0}=x\leq\theta_{2}^{*}$.
This establishes the contradiction. Hence, we can conclude that $E_{1}=\emptyset$
and $C_{1}=\mathcal{X}$, which implies that $\tau(E_{1})=\tau(\emptyset)$.
\\
(ii) Next, we prove the conditions under which $\{a_{1},a_{2}\}=\{\tau(E_{1}^{*}),\tau(\emptyset)\}$
is an MPE. Consider the following condition: 
\begin{align}
V_{2}(x;\tau(E_{1}^{*}),\tau(\emptyset))=\mathbb{E}^{x}\biggl[\int_{0}^{\tau_{1}^{*}}\pi(X_{t})e^{-rt}dt+w(X_{\tau_{1}^{*}}^{x})e^{-r\tau_{1}^{*}}\biggr]>l_{2} & \quad\textrm{ for all }x\in(\theta_{1}^{*},\theta_{2}^{*}]\:.\label{eq:Pf-pureMPE_Condition}
\end{align}

First, we prove that \eqref{eq:Pf-pureMPE_Condition} is a sufficient
condition for $\{a_{1},a_{2}\}=\{\tau(E_{1}^{*}),\tau(\emptyset)\}$
to be an MPE. Let $\tau(E_{2})$ be firm $2$'s best response to $\tau(E_{1}^{*})$,
i.e., $V_{W2}^{*}(x):=V_{2}(x;\tau(E_{1}^{*}),\tau(E_{2}))$ be the
corresponding payoff. We let $C_{2}=\mathcal{X}\backslash E_{2}$
denote the continuation region associated with the strategy $\tau(E_{2})$.

We now claim that $E_{2}\cap(\theta_{2}^{*},\infty)=\emptyset$. Towards
a contradiction, suppose not. Then we can pick some $x\in E_{2}\cap(\theta_{2}^{*},\infty)$,
which implies that $V_{W2}^{*}(x)=l_{2}$. However, because $\tau_{1}^{*}>\tau_{2}^{*}$
when $X_{0}=x$, Lemma \ref{lemm:Opt-Stop} implies that firm $2$
could obtain a \emph{strictly} higher payoff by exiting at $\tau_{2}^{*}>0$
instead, i.e.,
\begin{align*}
V_{W2}^{*}(x)\geq V_{2}(x;\tau(E_{1}^{*}),\tau(E_{2}^{*})) & =\mathbb{E}^{x}\biggl[\int_{0}^{\tau_{2}^{*}}\pi(X_{t})e^{-rt}dt+l_{2}e^{-r\tau_{2}^{*}}\biggr]>l_{2}\:,
\end{align*}
which is a contradiction. We next claim that $E_{2}\cap(\theta_{1}^{*},\theta_{2}^{*}]=\emptyset$.
Towards a contradiction, suppose not. Then we can pick $x\in E_{2}\cap(\theta_{1}^{*},\theta_{2}^{*}]$,
which implies that $V_{W2}^{*}(x)=l_{2}$. However, we have 
\begin{align*}
V_{W2}^{*}(x)\geq V_{2}(x;\tau(E_{1}^{*}),\tau(\emptyset)) & =\mathbb{E}^{x}\biggl[\int_{0}^{\tau_{1}^{*}}\pi(X_{t})e^{-rt}dt+w(X_{\tau_{1}^{*}}^{x})e^{-r\tau_{1}^{*}}\biggr]>l_{2}\:,
\end{align*}
where the last inequality follows from \eqref{eq:Pf-pureMPE_Condition}.
This establishes the contradiction. We further claim that $E_{2}\cap(-\infty,\theta_{1}^{*}]=\emptyset$.
If not, then there exists $x\in E_{2}\cap(-\infty,\theta_{1}^{*}]$,
which implies that both firms exit simultaneously when $X_{t}^{y}=x$,
and hence, $V_{W2}^{*}(x)=m_{2}(x)$. Because $\tau_{1}^{*}=0$ when
$X_{0}=x\leq\theta_{1}^{*}$, we have 
\begin{align*}
V_{W2}^{*}(x)\geq V_{2}(x;\tau(E_{1}^{*}),\tau(\emptyset)) & =\mathbb{E}^{x}\biggl[\int_{0}^{\tau_{1}^{*}}\pi(X_{t})e^{-rt}dt+w(X_{\tau_{1}^{*}}^{x})e^{-r\tau_{1}^{*}}\biggr]=w(x)>m_{2}(x)\:,
\end{align*}
which is a contradiction. Combining the three claims above, therefore,
we conclude that $E_{2}=\emptyset$, which implies that $C_{2}=\mathcal{X}$,
and hence, $\tau(E_{2})=\tau(\emptyset)$.

Second, define $\underline{w}:=\inf\{w(x):x\in\mathcal{X}\}$ and
$\beta_{W}(\theta):=[\underline{w}-R(\theta)]/\phi(\theta)$. Note
that $\beta_{W}(\theta)>\beta_{2}(\theta)$ for $\forall\theta\in\mathcal{X}$
because $\underline{w}>l_{2}$. Also, observe that for $\forall\theta<\theta_{2}^{*}$,
we have
\begin{align*}
\beta_{W}^{'}(\theta) & =\Bigl\{-R^{'}(\theta)\phi(\theta)-\phi^{'}(\theta)[\underline{w}-R(\theta)]\Bigr\}/\phi^{2}(\theta)\\
 & >\Bigl\{-R^{'}(\theta)\phi(\theta)-\phi^{'}(\theta)[l_{2}-R(\theta)]\Bigr\}/\phi^{2}(\theta)=\beta_{2}^{'}(\theta)>0
\end{align*}
where the first inequality follows because $\phi^{'}(\theta)<0$,
and the last inequality holds because $\beta_{2}^{'}(\theta)>0$ for
$\theta<\theta_{2}^{*}$ from Lemma \ref{lemm:Aux_1}. Next, pick
$\kappa_{\theta}>0$ such that 
\begin{align}
\beta_{W}(\theta_{2}^{*}-\kappa_{\theta}) & =\beta_{2}(\theta_{2}^{*})\:,\label{eq:Pf-pureMPE_Kappa}
\end{align}
where $\beta_{2}(\cdot)$ is defined in \eqref{eq:Beta_i}. If such
$\kappa_{\theta}$ exists, it must be unique because $\beta_{W}^{'}(\theta)>0$
for $\theta<\theta_{2}^{*}$. If there does not exist $\kappa_{\theta}$
which satisfies \ref{eq:Pf-pureMPE_Kappa}, then we let $\kappa_{\theta}=\infty$.

Finally, we show that \eqref{eq:Pf-pureMPE_Condition} is satisfied
$if$ $\theta_{2}^{*}-\theta_{1}^{*}<\kappa_{\theta}$, which will
complete the proof; this is because we can always find the unique
$\kappa_{l}>0$ for any given $\kappa_{\theta}>0$ such that $\theta_{2}^{*}-\theta_{1}^{*}<\kappa_{\theta}$
if and only if $l_{2}-l_{1}<\kappa_{l}$ from the fact that $\theta_{i}^{*}$
given in Lemma \ref{lemm:Opt-Stop} strictly increases in $l_{i}$.
Suppose now that $\theta_{2}^{*}-\theta_{1}^{*}<\kappa_{\theta}$,
i.e., $\theta_{1}^{*}>\theta_{2}^{*}-\kappa_{\theta}$. Note that
$\beta_{W}^{'}(\theta)>0$ for $\forall\theta<\theta_{2}^{*}$, and
recall that $\theta_{1}^{*}<\theta_{2}^{*}$. Therefore, $\beta_{W}(\theta_{1}^{*})>\beta_{W}(\theta_{2}^{*}-\kappa_{\theta})=\beta_{2}(\theta_{2}^{*})$
by \eqref{eq:Pf-pureMPE_Kappa}. Thus, for any $x\in(\theta_{1}^{*},\theta_{2}^{*}]$,
\begin{align*}
\mathbb{E}^{x}\biggl[\int_{0}^{\tau_{1}^{*}}\pi(X_{t}^{x})e^{-rt}dt+w(\theta_{1}^{*})e^{-r\tau_{1}^{*}}\biggr] & \geq\mathbb{E}^{x}\biggl[\int_{0}^{\tau_{1}^{*}}\pi(X_{t}^{x})e^{-rt}dt+\underline{w}e^{-r\tau_{1}^{*}}\biggr]=R(x)+\phi(x)\beta_{W}(\theta_{1}^{*})\\
 & >R(x)+\phi(x)\beta_{2}(\theta_{2}^{*})\geq R(x)+\phi(x)\beta_{2}(x)=l_{2}\:,
\end{align*}
where the first inequality holds from the definition of $\underline{w}$,
the first equality holds because $\mathbb{E}^{x}[e^{-r\tau_{1}^{*}}]=\phi(x)/\phi(\theta_{1}^{*})$
for $x>\theta_{1}^{*}$, the second inequality follows because $\beta_{W}(\theta_{1}^{*})>\beta_{2}(\theta_{2}^{*})$,
the last inequality holds because $\beta_{2}(\cdot)$ achieves its
maximum at $\theta_{2}^{*}$ by Lemma \ref{lemm:Aux_1}, and the last
equality follows by the definition of $\beta_{2}(\cdot)$. Hence,
\eqref{eq:Pf-pureMPE_Condition} is satisfied, which establishes the
desired result for $\kappa_{\theta}>0$.

\noindent (iii) Let $E_{i}\subset\mathbb{R}$ denote player $i$'s
exit region; \emph{i.e.}, it is player $i$'s strategy to exit whenever
$X_{t}$ enters $E_{i}$. Note that $\overline{E}_{1}\cap\overline{E}_{2}=\emptyset$
because none of the players have incentive to exit at the same time.
Without loss of generality, assume $\sup\left\{ E_{1}\right\} =\bar{x}_{1}>\sup\left\{ S_{2}\right\} =\bar{x}_{2}$.
Recall that we already established that $E_{i}\subseteq(-\infty,\theta_{i}^{*}]$,
so $\bar{x}_{1}\le\theta_{1}^{*}$. 

We now prove that $\bar{x}_{1}=\theta_{1}^{*}$ so that player 1's
equilibrium strategy is to exit at $\tau_{1}:=\inf\left\{ t>0\,:\,X_{t}\leq\theta_{1}^{*}\right\} $.
Towards a contradiction, suppose that $\bar{x}_{1}<\theta_{1}^{*}$
so that player 1's equilibrium strategy is to exit at $\bar{\tau}=\inf\{t\ge0:X_{t}=\bar{x}_{1}\}>\tau_{1}$.
Then player 1's payoff reduces to one in which player 1 exits at a
threshold $\bar{x}_{1}<\theta_{1}^{*}$ and player 2 never exits.
However, this contradicts Lemma 1, which asserts that the best response
of player 1 is to exit at threshold $\theta_{1}^{*}$. We conclude
that $\bar{x}_{1}=\theta_{1}^{*}$ if $\bar{x}_{1}>\bar{x}_{2}$,
which proves the statement.\eproof

\textbf{Proof of Lemma \ref{lemma:MixedMPE2}}: Suppose that $\left\{ a_{1},a_{2}\right\} =\left\{ \left(E_{1},\lambda_{1}\right),\left(E_{2},\lambda_{2}\right)\right\} $
is a mixed-strategy MPE, where $E_{i}$ is firm $i$'s exit region
(\emph{i.e.}, it exits with probability one whenever $X_{t}\in E_{i}$),
and $\lambda_{i}(x)$ is its hazard rate of exit function. For the
sake of generality, we do not assume that $E_{i}$ is a closed set
at the outset in this proof. Instead, we prove that $E_{i}$ can be
replaced by its closure without affecting the payoff functions (see
the remark followed by Claim 2 below).

Let $\Gamma_{i}$ denote the support of firm $i$'s mixed strategy
as defined in (\ref{eq:supp-lambda}). Without loss of generality,
we can assume that $E_{i}\cap\Gamma_{i}=\emptyset$ because if firm
$i$ exits at some $x\in\mathcal{X}$, it does so either with probability
$1$ (if $x\in E_{i}$) or at the hazard rate $\lambda_{i}(x)$ (if
$x\in\Gamma_{i}$) based on our definition of Markov strategy.

We first show that $E_{i}\cap E_{-i}=\emptyset$ and $\Gamma_{i}\cap E_{-i}=\emptyset$.
For the purpose of contradiction, suppose that there exists some $y\in E_{i}\cap E_{-i}$.
Then firm $i$'s payoff at $y$ is $\frac{1}{2}(l_{i}+w(y))$ because
both firms attempt to exit simultaneously. However, if firm $i$ deviates
from its strategy and chooses not to exit, its payoff at $y$ would
be $w(y)>\frac{1}{2}(l_{i}+w(y))$, so its payoff is improved. This
contradicts the assumption that $\left\{ a_{1},a_{2}\right\} $ is
an equilibrium, and hence, we conclude $E_{i}\cap E_{-i}=\emptyset$.

Next, suppose that there exists $y\in\Gamma_{i}\cap E_{-i}$. Note
that by the definition of a mixed strategy equilibrium, one of firm
$i$'s best responses to $a_{-i}$ is to attempt to exit immediately
whenever the current value of $X$ is within $\Gamma_{i}$. Since
$y$ also belongs to $E_{-i}$, firm $i$'s payoff at $y$ is $\frac{1}{2}(l_{i}+w(y))$,
which is strictly less than the payoff from not exiting at $y$, i.e.,
$w(y)$. This implies that firm $i$ can improve its payoff by deviating
from its candidate equilibrium strategy, which again contradicts the
assumption that $\left\{ a_{1},a_{2}\right\} $ is an equilibrium.
Thus, we conclude $\Gamma_{i}\cap E_{-i}=\emptyset$. 

Now, we prove two claims that will be used for the proof of the lemma.
Below we let $\text{cl}(S)$ denote the closure of the set $S$ and
$\partial S$ denote the boundary of the set $S$.

\medskip{}

\textbf{Claim 1}: $\Gamma_{i}=\Gamma_{-i}$ and $\sup\Gamma_{i}\le x_{ci}$.
Furthermore, $\lambda_{i}(x)=\frac{rl_{-i}-\pi(x)}{w(x)-l_{-i}}$
for any $x\in\Gamma_{i}$.

Proof: Recall that $\Gamma_{i}$ can be expressed as a union of open
intervals. Towards a contradiction, suppose that for some $i$, there
exists an open interval $G$ such that $G\subseteq\Gamma_{i}$ and
$G\cap\Gamma_{-i}=\emptyset$. We let $D_{i}:=\{x\in\mathcal{X}:\pi(x)>rl_{i}\}=(x_{ci},\beta)$
denote the set of states at which firm $i$'s net-present-value of
its flow payoff from remaining in the market exceeds its outside option.
Without loss of generality, we can always appropriately choose $G$
such that either $G\subset D_{i}$ or $G\cap D_{i}=\emptyset$ holds.

We first consider the case $G\subset D_{i}$. Assume $X_{0}\in G$
and define $\tau_{G}:=\inf\left\{ t\,:\,X_{t}\notin G\right\} $.
Then $\tau_{G}>0$ a.s because $G$ is an open set. Now, fix some
stopping time $\tau\in(0,\tau_{G})$. Then because $X_{0},X_{\tau}\in G$
and $G$ is a subset of $\Gamma_{i}$, firm $i$'s payoff from an
immediate exit should be the same as that from an exit at $\tau$
by the definition of a mixed strategy equilibrium. However, because
$G\cap\Gamma_{-i}=\emptyset$ by assumption, firm $i$'s expected
payoff if it exits at $\tau$ and its rival's strategy is $a_{-i}$
is equal to
\begin{equation}
l_{i}+\mathbb{E}\biggl[\int_{0}^{\tau}e^{-rt}\left(\pi(X_{t})-rl_{i}\right)dt\biggr]>l_{i}\,,\label{eq:Vi-tau}
\end{equation}
where the inequality follows from the fact that $\pi(X_{t})>rl_{i}$
for all $t\in[0,\tau)$. Therefore, firm $i$'s payoff is strictly
higher if it exits at $\tau>0$ than if it exits immediately, which
contradicts the assumption that $\{a_{1},a_{2}\}$ is a mixed strategy
equilibrium.

Next, we consider the case $G\cap D_{i}=\emptyset$. We repeat the
same argument as the case of $G\subset D_{i}$ with $\tau\in(0,\tau_{G})$.
Then firm $i$'s expected payoff if it exits at $\tau$ and its rival's
strategy is $a_{-i}$ is equal to
\[
l_{i}+\mathbb{E}\biggl[\int_{0}^{\tau}e^{-rt}\left(\pi(X_{t})-rl_{i}\right)dt\biggr]<l_{i}\:,
\]
where the inequality is due to $\pi(X_{t})<rl_{i}$ for all $t\in[0,\tau)$.
Therefore, firm $i$ can obtain a strictly higher payoff from an immediate
exit than an exit at $\tau$. This contradicts the assumption that
$\{a_{1},a_{2}\}$ is a mixed strategy equilibrium. We conclude that
$\Gamma_{i}=\Gamma_{-i}$ must hold.

Now we derive the form of the equilibrium rate of exit. For notational
convenience, we let $V(x)$ denote the payoff function $V_{i}(x;a_{i},a_{-i})$.
Recall that $V(x)=l_{i}$ for $x\in\Gamma:=\Gamma_{i}=\Gamma_{-i}$
because an immediate exit within $\Gamma$ is each firm's best response.
Let $G$ be an open neighborhood of $x$ such that $G\subsetneq\Gamma$,
and let $\tau_{G}:=\inf\left\{ t\ge0\,:\,X_{t}\notin G\right\} $.
By the definition of the mixed strategy region $\Gamma$, firm $i$'s
payoff associated with exit at time $\tau_{G}$ is still $V(x)=l_{i}$
because an exit at $\tau_{G}$ is also one of the best responses of
firm $i$. Hence, one of the possible expressions of $V(x)$ is given
by the following:
\begin{equation}
V(x)=\mathbb{E}^{x}[\int_{0}^{\tau_{G}}e^{-\Lambda_{t}^{(-i)}}(\pi(X_{t})+\lambda_{-i}(X_{t})w(X_{t}))dt+e^{-\Lambda_{\tau_{G}}^{(-i)}}V(X_{\tau_{G}})]\:,\label{eq:V-alternative}
\end{equation}
where $\Lambda_{t}^{(-i)}:=\int_{0}^{t}(\lambda_{-i}(X_{s})+r)ds$.
Intuitively, the discount factor $e^{-rt}$ is replaced by $e^{-\Lambda_{t}^{(-i)}}$
due to the exit rate $\lambda_{-i}(\cdot)$ of firm $-i$; the additional
flow profit $\lambda_{-i}(X_{t})w(X_{t})$ is because of the exit
probability of firm $-i$; the reward $l_{i}$ from exit at $\tau_{G}$
is replaced by $V(X_{\tau_{G}})$ because $V(X_{\tau_{G}})=l_{i}$.

To derive the form of $\lambda_{-i}(\cdot)$, we need to derive the
Hamilton-Jacobi-Bellman (HJB) equation:
\begin{equation}
[\frac{\sigma(x)^{2}}{2}\frac{d^{2}}{dx^{2}}+\mu(x)\frac{d}{dx}-r-\lambda_{-i}(x)]V(x)+\lambda_{-i}(x)w(x)+\pi(x)=0\:.\label{eq:HJB-eq}
\end{equation}
 For this purpose, we define a function $f(\Lambda,x):=e^{-\Lambda}V(x)$.
The new variable $\Lambda$ represents the integral $\Lambda_{t}^{(-i)}$
which satisfies the differential equation $d\Lambda_{t}^{(-i)}=(r+\lambda_{-i}(X_{t}))dt$
by its definition. Note that $f(\Lambda,x)$ is continuously differentiable
with respect to $\Lambda$ and twice continuously differentiable with
respect to $x$ within $\Gamma$ because $V(x)=l_{i}$ for $x\in\Gamma$.
Thus, based on the stochastic differential equation that the bivariate
process $(\Lambda_{t}^{(-i)},X_{t})$ follows, the following characteristic
operator is well-defined for the function $f$:
\[
\mathcal{L}:=(\lambda_{-i}(x)+r)\frac{\partial}{\partial\Lambda}+\frac{\sigma(x)^{2}}{2}\frac{\partial^{2}}{\partial x^{2}}+\mu(x)\frac{\partial}{\partial x}\;.
\]

From the expression (\ref{eq:V-alternative}), we have the following:
\begin{equation}
\mathbb{E}^{x}[f(\Lambda_{\tau_{G}}^{(-i)},X_{\tau_{G}})]-f(0,x)=-\mathbb{E}^{x}[\int_{0}^{\tau_{G}}e^{-\Lambda_{t}^{(-i)}}(\pi(X_{t})+\lambda_{-i}(X_{t})w(X_{t}))dt]\:.\label{eq:V-HJB}
\end{equation}
Recall that $G$ is any arbitrary open neighborhood of $x$. From
equation (7.5.1) of \citet{Oksendal03}, 
\begin{align}
\lim_{G\downarrow x}\frac{\mathbb{E}^{x}[f(\Lambda_{\tau_{G}}^{(-i)},X_{\tau_{G}})]-f(0,x)}{\mathbb{E}^{x}[\tau_{G}]} & =\mathcal{L}f(0,x)=[\frac{\sigma(x)^{2}}{2}\frac{d^{2}}{dx^{2}}+\mu(x)\frac{d}{dx}-r-\lambda_{-i}(x)]V(x)\:.\label{eq:HJB-LHS}
\end{align}
Here the limit $G\downarrow x$ represents a series of open sets decreasing
to the point $x$. Furthermore, 
\begin{equation}
\lim_{G\downarrow x}\frac{-\mathbb{E}^{x}[\int_{0}^{\tau_{G}}e^{-\Lambda_{t}^{(-i)}}(\pi(X_{t})+\lambda_{-i}(X_{t})w(X_{t}))dt]}{\mathbb{E}^{x}[\tau_{G}]}=-\pi(x)-\lambda_{-i}(x)w(x)\:.\label{eq:HJB-RHS}
\end{equation}
Thus, from (\ref{eq:V-HJB}), (\ref{eq:HJB-LHS}), and (\ref{eq:HJB-RHS}),
we obtain (\ref{eq:HJB-eq}).

Because $V(\cdot)=l_{i}$, the HJB equation (\ref{eq:HJB-eq}) turns
into $[-r-\lambda_{-i}(x)]l_{i}+\lambda_{-i}(x)w(x)+\pi(x)=0$, which
leads to 
\begin{equation}
\lambda_{-i}(x)=\frac{rl_{i}-\pi(x)}{w(x)-l_{i}}\:.\label{eq:lambda-i}
\end{equation}
Finally, we impose the condition that $\lambda_{-i}(\cdot)$ must
be positive within $\Gamma_{-i}$, so we obtain $rl_{i}>\pi(x)$ because
$w(\cdot)>l_{i}$. Therefore, $\sup\Gamma_{i}\le x_{ci}$ and $\lambda_{i}(x)=\frac{rl_{-i}-\pi(x)}{w(x)-l_{-i}}$
for $x\in\Gamma_{i}$. $\square$

\medskip{}

\textbf{Claim 2}: $\partial\Gamma_{i}\cap\partial E_{-i}=\emptyset$
and $\partial E_{i}\cap\partial E_{-i}=\emptyset$.

Proof: Suppose that there exists $y\in\partial\Gamma_{i}\cap\partial E_{-i}$.
Since $\Gamma_{i}$ is a union of disjoint open intervals, there is
a subset of $\Gamma_{i}$ of the form $[z,y)$ or $(y,z]$. Without
loss of generality, we assume that $[z,y)\subset\Gamma_{i}$ for some
$z<y$.

From Claim 1, $[z,y)$ is also a subset of $\Gamma_{-i}$, within
which firm $-i$ exits with a rate $\lambda_{-i}(x)$ given by (\ref{eq:lambda-i}).
Suppose that, within the interval of $[z,y)$, firm $i$ employs an
alternative strategy $a_{i}'$ of exit at the escape time $\tau=\inf\{t\ge0:X_{t}\not\in(z,y')\}$
of $(z,y')$ for some $y'>y$. Then we have $V_{i}(z;a_{i}',a_{-i})=l_{i}$
by the definition of the strategy $a_{i}'$. Furthermore, $V_{i}(y;a_{i}',a_{-i})=w(y)$
because firm $-i$ exits at the hitting time of $y\in\partial E_{-i}$
due to the property of the regular diffusion process $X$. Along with
the two boundary conditions, because $\lambda_{-i}(\cdot)$ is a continuous
function within $\Gamma_{-i}$, $V_{i}(x;a_{i}',a_{-i})$ for $x\in(z,y)$
can be obtained as the solution to a boundary value problem, satisfying
the following differential equation (Section 9.1, \citealt{Oksendal03}):
\[
[\frac{\sigma(x)^{2}}{2}\frac{d^{2}}{dx^{2}}+\mu(x)\frac{d}{dx}-r-\lambda_{-i}(x)]V_{i}(x;a_{i}',a_{-i})+\lambda_{-i}(x)w(x)+\pi(x)=0\:.
\]
It is well-known that there exists a unique solution $V_{i}(x;a_{i}',a_{-i})$
to this differential equation within $(z,y)$ (see, for example, Appendix
of \citealt{Lon2011}), and it follows that $V_{i}(x;a_{i}',a_{-i})$
is continuous within $[z,y]$. From the boundary condition $V_{i}(y;a_{i}',a_{-i})=w(y)>l_{i}$
and the continuity of $V_{i}(\cdot;a_{i}',a_{-i})$, there must exist
an interval $I\subset[z,y]$ such that $V_{i}(x;a_{i}',a_{-i})>l_{i}$
for all $x\in I$. This implies that firm $i$ can improve its payoff
within $I$ by adopting an alternative strategy $a_{i}'$, which contradicts
the assumption that $\{a_{i},a_{-i}\}$ is an equilibrium. We conclude
that $\partial\Gamma_{i}\cap\partial E_{-i}=\emptyset$.

Next, suppose that there exists $y\in\partial E_{i}\cap\partial E_{-i}$.
Because $\Gamma_{-i}\cap E_{i}=\emptyset$ and $\partial\Gamma_{-i}\cap\partial E_{i}=\emptyset$,
there exists an open neighborhood $G$ of $y$ such that $G\cap\Gamma_{-i}=\emptyset$.
Within $G$, we consider an alternative strategy $a_{i}'$ of exit
at the escape time $\tau_{G}:=\inf\left\{ t\ge0\,:\,X_{t}\notin G\right\} $
of $G$. Then for any $x\in G\backslash\text{cl}(E_{-i})$, where
$G\backslash\text{cl}(E_{-i})$ is an open set, firm $i$'s payoff
function $V_{i}(x;a_{i}',a_{-i})$ within $G\backslash\text{cl}(E_{-i})$
can be obtained as the solution to a boundary value problem (Section
9.1, \citealt{Oksendal03}), satisfying the following differential
equation
\[
[\frac{\sigma(x)^{2}}{2}\frac{d^{2}}{dx^{2}}+\mu(x)\frac{d}{dx}-r]V_{i}(x;a_{i}',a_{-i})+\pi(x)=0\:,
\]
along with the boundary conditions (i) $V_{i}(x;a_{i}',a_{-i})=l_{i}$
at $x\in\partial G\backslash\text{cl}(E_{-i})$, (ii) $V_{i}(x;a_{i}',a_{-i})=(l_{i}+w(x))/2$
at $x\in\partial G\cap E_{-i}$, and (iii) $V_{i}(x;a_{i}',a_{-i})=w(x)$
for all $x\in G\cap\text{cl}(E_{-i})$ because firm $-i$ exits at
the hitting time of $\text{cl}(E_{-i})$ by the diffusive property
of $X$. By the solution property of the boundary value problem, $V_{i}(x;a_{i}',a_{-i})$
is continuous within the set $\text{cl}(G)$. Recall that $y\in G\cap\partial E_{-i}$,
which implies that $V_{i}(y;a_{i}',a_{-i})=w(y)>l_{i}$ from case
(iii) above, and that $y\in\partial E_{i}$; hence, we can choose
$x\in E_{i}$ arbitrarily close to $y$ such that $V_{i}(x;a_{i}',a_{-i})>l_{i}$
because of the continuity of $V_{i}(\cdot;a_{i}',a_{-i})$. Therefore,
because $V_{i}(x;a_{i},a_{-i})=l_{i}<V_{i}(x;a_{i}',a_{-i})$ for
some $x\in E_{i}$, firm $i$ can improve the payoff by deviating
from $a_{i}$ and adopting an alternative strategy $a_{i}'$, which
contradicts the assumption that $\{a_{i},a_{-i}\}$ is an equilibrium.
We conclude that $\partial E_{i}\cap\partial E_{-i}=\emptyset$. $\square$

\medskip{}

\emph{Remark}: Because we show $\partial E_{i}\cap\partial E_{-i}=\emptyset$
based on the argument that holds for any exit probability $p_{i}\in(0,1]$
at $E_{i}$, $i\in\{1,2\}$, we do not need to consider the dependence
of the payoffs on $p_{i}$'s at a point $y\in\partial E_{i}\cap\partial E_{-i}$
in our equilibrium analysis. Also, by Claim 2 and the fact that $E_{i}\cap E_{-i}=\emptyset$
and $\Gamma_{i}\cap E_{-i}=\emptyset$, we established that both $E_{i}$
and $\Gamma_{i}$ must be separated from $E_{-i}$ by closed neighborhoods
(\emph{i.e.}, there must exist a non-empty interval between the two
sets). Therefore, without loss of generality, we can choose $E_{i}$
and $E_{-i}$ to be closed sets without affecting the payoff functions
because the hitting time of $E_{i}$ is identical to the hitting time
of $\text{cl}(E_{i})$ due to the diffusive property of $X$.

With the two claims established above, we now prove the statement
of the lemma:

\medskip{}

\emph{Step 1.---}For each $i$, let $C_{i}$ denote the set of states
at which firm $i$ does not exit, that is, $C_{i}:=(\alpha,\beta)\backslash(\Gamma_{i}\cup E_{i})$.
Moreover, let $F_{i}:=(\theta_{i}^{*},\beta)$ denote the set of states
at which firm $i$ would prefer to remain in the market if it expects
its rival to never exit. Recall from Lemma \ref{lemm:Opt-Stop} that
$D_{i}\subset F_{i}$. We will show that $F_{i}\subset C_{i}$, or
equivalently, $(\Gamma_{i}\cup E_{i})\cap F_{i}=\emptyset$, that
is, $F_{i}$ is a subset of the continuation region for firm $i$.
Towards a contradiction, suppose there exists some $x\in(\Gamma_{i}\cup E_{i})\cap F_{i}$.
Fix an $i$, and define the strategies $\widetilde{a}_{i}:=\left((\alpha,\theta_{i}^{*}],\mathbf{0}\right)$
and $\widetilde{a}_{-i}:=\left(\emptyset,\mathbf{0}\right)$, that
is, a strategy profile in which firm $i$ exits (with probability
1) whenever $X_{t}\in(\alpha,\theta_{i}^{*}]$ and its rival never
exits. Then
\[
V_{i}(x;a_{i},a_{-i})=l_{i}<V_{i}(x;\widetilde{a}_{i},\widetilde{a}_{-i})\leq V_{i}(x;\widetilde{a}_{i},a_{-i}).
\]
The equality follows from the assumption that $x\in\Gamma_{i}\cup E_{i}$,
which also implies $x\notin E_{-i}$. The first inequality follows
because $x>\theta_{i}^{*}$ (since $x\in F_{i}$ by assumption), so
by Lemma \ref{lemm:Opt-Stop}, firm $i$ can obtain a strictly higher
payoff by exiting at $\theta_{i}^{*}$ if its rival never exits. The
last inequality follows because firm $i$ is better off if its rival
exits in finite time with positive probability compared to the case
in which it never exits. To elaborate, the assumption that $w(x)>\mathbb{E}[\int_{0}^{t}e^{-rs}\pi(X_{s})ds+e^{-rt}w(X_{t})|X_{0}=x]$
for all $x\in\mathcal{X}$ and $t$ implies that the payoff process
for the winner $W_{i}(t):=\int_{0}^{t}\pi(X_{s})e^{-rs}ds+e^{-rt}w(X_{t})$
is a supermartingale. Letting $\tau_{\theta_{i}^{*}}:=\inf\{t:X_{t}\le\theta_{i}^{*}\}$
denote the first hitting time of the set $(\alpha,\theta_{i}^{*}]$,
the supermartingale property of $W_{i}(\cdot)$ implies that $W_{i}(t)\geq\mathbb{E}[W_{i}(\tau_{\theta_{i}^{*}})|\mathcal{F}_{t}]$
for any $t<\tau_{\theta_{i}^{*}}$, that is, firm $i$ is better off
(in expectation) becoming the winner at any $t<\tau_{\theta_{i}^{*}}$
compared to becoming the winner at $t=\tau_{\theta_{i}^{*}}$, which
is in turn strictly better than becoming the loser at $t=\tau_{\theta_{i}^{*}}$
because $w(X_{\tau_{\theta_{i}^{*}}})>l_{i}$. This implies that firm
$i$'s expected payoff $V_{i}(x;\widetilde{a}_{i},\widetilde{a}_{-i})\leq V_{i}(x;\widetilde{a}_{i},a_{-i})$
for any strategy $a_{-i}$. We have thus established that $\Gamma_{i}\cup E_{i}$
does not intersect with $D_{i}$ or $F_{i}$.

\medskip{}

\emph{Step 2.---}Next, we prove that $E_{1}=E_{2}=\emptyset$. Towards
a contradiction, suppose that $E_{i}\ne\emptyset$. Then because each
$E_{i}$, $i\in\{1,2\}$, is a closed set, so must be their union
$E_{1}\cup E_{2}$, which implies that its complement $\mathcal{X}\backslash(E_{1}\cup E_{2})$
is an open set in the state space $\mathcal{X}=(\alpha,\beta)\subseteq\mathbb{R}$.
Hence, $\mathcal{X}\backslash(E_{1}\cup E_{2})$ must be a union of
open intervals because any open set on the real line can be expressed
as a union of open intervals. Now, we can always find a subinterval
$(f,g)\subseteq\Gamma$ such that $(f,g)$ is a \emph{proper} subset
of an open interval contained in $\mathcal{X}\backslash(E_{1}\cup E_{2})$;
in other words, there exists an interval $(c,d)\subseteq\mathcal{X}\backslash(E_{1}\cup E_{2})$
such that $c$ or $d\in E_{1}\cup E_{2}$ and $c<f<g<d$. This is
because (1) $\Gamma_{i}$ is always a union of intervals by assumption,
(2) $\Gamma=\Gamma_{1}=\Gamma_{2}$ by Claim 1, (3) $\Gamma_{i}$
is separated from $E_{-i}$ by closed neighborhoods by Claim 2.

Suppose first that $d\in E_{i}$ for some $i\in\{1,2\}$. Without
loss of generality, we can always choose $(f,g)$ in a way that $(g,d)$
is a subset of the continuation region for firm $i$; this is because
$(f,g)\subseteq(c,d)\subseteq\mathcal{X}\backslash(E_{1}\cup E_{2})$
by construction and $\Gamma$ is separated from $E_{i}$ by closed
neighborhoods (so any components of $\Gamma$ contained in $(c,d)$
must be strictly below $d$). Also, because we proved $\Gamma_{i}\cup E_{i}$,
$i\in\{1,2\}$, does not intersect $F_{i}=(\theta_{i}^{*},\beta)$
in Step 1, it follows that $g\leq\min\{\theta_{1}^{*},\theta_{2}^{*}\}$
(since $\Gamma_{1}=\Gamma_{2}$) and $d\leq\theta_{i}^{*}$. Moreover,
because $(f,g)\subseteq\Gamma$ and $d\in E_{i}$, we have $V_{i}(d;a_{1},a_{2})=V_{i}(g;a_{1},a_{2})=l_{i}$.

\emph{Remark}: The boundary conditions $V_{i}(d;a_{1},a_{2})=V_{i}(g;a_{1},a_{2})=l_{i}$
can be derived as follows. First, $V_{i}(d;a_{1},a_{2})=l_{i}$ follows
from the fact that firm $i$ exits at the hitting time of $d\in E_{i}$
under the strategy $a_{i}$. Second, $V_{i}(g;a_{1},a_{2})=l_{i}$
holds for the following reason. If $(a_{1},a_{2})$ is a mixed strategy
MPE, by the definition of a mixed strategy equilibrium, firm $i$'s
payoff function must remain unchanged even if firm $i$ employs an
alternative strategy of pure-strategy exit within $\Gamma$. In this
case, firm $i$ exits at the hitting time of $g$, which is the boundary
point of $\Gamma$. (Because $X$ is a diffusion process, the hitting
time of $\Gamma$ is identical to the hitting time of the closure
of $\Gamma$.) It thus follows that firm $i$'s payoff at $g$ is
$l_{i}$. Alternatively, one can invoke the continuity of $V_{i}(x;a_{1},a_{2})$
at $x=g$ to arrive at the boundary condition $V_{i}(g;a_{1},a_{2})=l_{i}$
because $V_{i}(x;a_{1},a_{2})=l_{i}$ for all $x\in\Gamma$.

Now, fix some $X_{0}=x\in(g,d)$, and assume that firm $i$ exits
at the hitting time $\tau_{b}=\inf\left\{ t\,:\,X_{t}\in\{g,d\}\right\} $;
\emph{i.e.}, the first time that $X_{t}$ hits $\{g,d\}$. Define
$a_{i}'=(\{g,d\},\mathbf{0})$. Then we have $V_{i}(x;a_{i},a_{-i})=V_{i}(x;a_{i}',a_{-i})$.
Because $\pi(X_{t})<rl_{i}$ for all $t<\tau_{b}$ (by $g\leq\theta_{i}^{*}$
and $d\leq\theta_{i}^{*}$), we can get the inequality similar to
(\ref{eq:Vi-tau}) to conclude that $V_{i}(x;a_{i},a_{-i})<l_{i}$,
which contradicts the assumption that $a_{i}$ is a best response
to $a_{-i}$.

Suppose next that $c\in E_{i}$ for some $i\in\{1,2\}$. Then we can
similarly proceed with firm $i$'s continuation region $(c,f)$ and
leads to a contradiction to the assumption that $a_{i}$ is a best
response to $a_{-i}$. Therefore, we conclude that $E_{1}=E_{2}=\emptyset$.\medskip{}

\emph{Step 3.---} Finally, we prove that $\Gamma=(\alpha,\theta_{i}^{*})$.
Towards doing so, we will first show that $\Gamma=(\alpha,\theta)$
for some $\theta$; note that, if this is the case, we must have $\theta\leq\min\{\theta_{1}^{*},\theta_{2}^{*}\}$
because $\Gamma_{i}\cap F_{i}=\emptyset$ and $\Gamma=\Gamma_{1}=\Gamma_{2}$
by Step 1 and Claim 1 respectively. Let $J:=\mathcal{X}\backslash\Gamma.$
We now consider the following two cases for proof by contradiction.

(i) First, suppose that there exists an interval $(c,d)$ such that
$\alpha<c<d<\min\{\theta_{1}^{*},\theta_{2}^{*}\}$ and $(c,d)\subseteq J$
with $c,d\in\text{cl}(\Gamma)$, i.e., $(c,d)$ belongs to both firms'
continuation region, and yet, it borders with $\Gamma$ at both ends.
Let $X_{0}=x\in(c,d)$ and consider the hitting time $\tau_{(c,d)}:=\inf\{t:X_{t}\notin(c,d)\}$,
\emph{i.e., }the first time that $X_{t}$ hits the set $\{c,d\}$.
Similarly as in Step 2, by the definition of a mixed strategy equilibrium,
firm $i$'s payoff function remains unchanged even if he employs an
alternative pure strategy of exit at $\tau_{(c,d)}$. In other words,
if we define $a_{i}'=(\{c,d\},\mathbf{0})$ as the said alternative
strategy, then we have $V_{i}(x;a_{i},a_{-i})=V_{i}(x;a_{i}',a_{-i})$,
which satisfies the following inequality:
\begin{equation}
V_{i}(x;a_{i}',a_{-i})=l_{i}+\mathbb{E}\left[\int_{0}^{\tau_{(c,d)}}[\pi(X_{t})-rl_{i}]e^{-rt}dt\right]<l_{i}\,,\label{eq:Vi-Step4}
\end{equation}
where the inequality follows because $\pi(X_{s})<rl_{i}$ for all
$s\leq\tau_{(c,d)}$. Therefore firm $i$ is strictly better off employing
a strategy of an immediate exit--a contradiction.

(ii) Next, suppose that there exists an interval $(\alpha,c)$ such
that $\alpha<c<\min\{\theta_{1}^{*},\theta_{2}^{*}\}$ and $(\alpha,c)\subseteq J$
with $c\in\text{cl}(\Gamma)$, i.e., $(\alpha,c)$ belongs to both
firms' continuation region yet it borders with $\Gamma$ at the point
$c$. Let $X_{0}=x\in(\alpha,c)$ and consider the hitting time $\tau_{c}:=\inf\{t:X_{t}\geq c\}$,
\emph{i.e., }the first time that $X_{t}$ hits the point $c$. Then
similarly as in the case of $(c,d)\subseteq J$ above, $V_{i}(x;\left(\{c\},\mathbf{0}\right),a_{-i})$
must be equal to firm $i$'s equilibrium payoff $V_{i}(x;a_{i},a_{-i})$,
but we will get the inequality similar to \eqref{eq:Vi-Step4} for
$V_{i}(x;\left(\{c\},\mathbf{0}\right),a_{-i})$, which leads to a
contradiction; here, although $\mathbb{P}(\tau_{c}=\infty)>0$ is
possible, $V_{i}(x;\left(\{c\},\mathbf{0}\right),a_{-i})$ is still
well-defined because we assume in Section \ref{sec:Model} that $\pi(\cdot)$
satisfies the absolute integrability condition.

(iii) From (i) and (ii), we conclude that $\Gamma=(a,\theta)$ for
some $\theta$.

Finally, we show that $\theta=\theta_{i}^{*}$ for each $i$. Recall
from Step 1 that $\Gamma_{i}\cup E_{i}$ does not intersect with $F_{i}=(\theta_{i}^{*},\beta)$
for any $i$, from Claim 1 that $\Gamma_{1}=\Gamma_{2}=\Gamma$, and
from Step 2 that $E_{1}=E_{2}=\emptyset$. Therefore, $\Gamma\subseteq(\alpha,\theta_{i}^{*}]$,
and so $\theta\leq\theta_{i}^{*}$ must hold for each $i$. Towards
a contradiction, suppose that $\theta<\theta_{i}^{*}$ and fix some
$x\in(\theta,\theta_{i}^{*})$. Letting $\tau_{\theta}=\inf\left\{ t\,:\,X_{t}\leq\theta\right\} $,
notice that
\[
V_{i}(x;a_{1},a_{2})=\mathbb{E}\left[\int_{0}^{\tau_{\theta}}e^{-rt}\pi(X_{t})dt+e^{-r\tau_{\theta}}l_{i}\,|\,X_{0}=x\right]=l_{i}+\mathbb{E}\left[\int_{0}^{\tau_{\theta}}e^{-rt}\left(\pi(X_{t})-rl_{i}\right)dt\,|\,X_{0}=x\right]<l_{i}\,.
\]
The first equality follows from the fact that $V_{i}(\theta;a_{1},a_{2})=l_{i}$
since firm $i$ is indifferent between exiting and remaining in the
market when $x=\theta$, the second equality follows by manipulating
terms, and the inequality follows from the fact that $\pi(X_{t})<rl_{i}$
for all $t<\tau_{\theta}$ (recall from Step 1 that $D_{i}\subset F_{i}$).
Therefore, firm $i$ is strictly better off exiting immediately, contradicting
the premise that $\{a_{1},a_{2}\}$ is an MPE. Hence we conclude that
$\theta=\theta_{i}^{*}$ for each $i$.\eproof

\noindent \textbf{Proof of Theorem} \textbf{\ref{thm:Stoch-MixedMPE}}:
By noting that $\theta_{1}^{*}=\theta_{2}^{*}$ if and only if $l_{1}=l_{2}$,
it follows immediately from Lemma \ref{lemma:MixedMPE2} that if $l_{1}<l_{2}$,
then the game does not admit any mixed strategy MPE. \eproof

\newpage{}

\section{Online Appendix: Structural Stability \label{sec:Structural-Stability}}

The goal of this section is to investigate the robustness of Theorem
\ref{thm:Stoch-MixedMPE} to our assumptions regarding the firms'
payoffs. In particular, we argue that Theorem \ref{thm:Stoch-MixedMPE}
continues to hold even if the firms have heterogeneous discount rates
($r_{1}\neq r_{2}$), heterogeneous flow profits while they remain
in the market ($\pi_{1}(\cdot)\neq\pi_{2}(\cdot)$), heterogeneous
winner payoffs ($w_{1}(\cdot)\neq w_{2}(\cdot)$), and the loser's
payoff is state-dependent (\emph{i.e.}, if firm $i$ exits at $t$,
then it obtains payoff $l_{i}(X_{t})$).

To analyze this model, in addition to the assumptions at the end of
the model description in Section \ref{sec:Model}, we make the following
assumptions:\begin{condition}\label{Ass-Str-Stability} Assume that
for each $i\in\{1,2\}$,

\noindent (i) $w_{i}(x)>l_{i}(x)$ for all $x\in\mathcal{X}$,

\noindent (ii) $l_{i}(\cdot)$ is twice continuously differentiable
on $\mathcal{X}$,

\noindent (iii) $\pi_{i}(x)+\mathcal{A}_{i}l_{i}(x)$ is increasing
in $x$, and

\noindent (iv) $\lim_{x\downarrow a}\pi_{i}(x)+\mathcal{A}_{i}l_{i}(x)<0$
and $\lim_{x\uparrow b}\pi_{i}(x)+\mathcal{A}_{i}l_{i}(x)>0$.\end{condition}

Part (i) ensures that the winner's payoff is always greater than that
of the loser, and it is analogous to the assumption $w(\cdot)>l_{2}$
we made in Section \ref{sec:Model}. Part (ii) implies that we can
apply the infinitesimal generator $\mathcal{A}_{i}$ to $l_{i}(\cdot)$.
Parts (iii) and (iv) guarantee that there exists a unique threshold
$\theta_{i}^{*}$ such that firm $i$'s best response if its rival
never exits (\emph{i.e.}, if $a_{-i}=\left\{ \emptyset,\emptyset,\mathbf{0}\right\} $)
is to exit at the first time $\tau_{i}=\inf\{t\geq0:X_{t}\leq\theta_{i}^{*}\}$
(see Theorem 6 (B) in \citet{Alvarez2001} for details).\footnote{As an example, if the state $X$ is a linear diffusion (\emph{i.e.},
$\mu(x)\equiv\mu<0$ and $\sigma(x)\equiv\sigma>0$ in (\ref{eq:State})),
$\pi_{i}(x)=A_{i}x+B_{i}$ and $l_{i}(x)=C_{i}x+D_{i}$, then it is
easy to verify that Condition \ref{Ass-Str-Stability}(ii)-(iv) are
satisfied as long as $A_{i}>r_{i}C_{i}$. If $X$ is a geometric Brownian
motion (\emph{i.e.}, $\mu(x)\equiv\mu x$ and $\sigma(x)\equiv\sigma x$
for some $\mu<0$ and $\sigma>0$), then Condition \ref{Ass-Str-Stability}(ii)-(iv)
are satisfied for the above choice of $\pi_{i}(\cdot)$ and $l_{i}(\cdot)$
as long as $A_{i}>(r_{i}+\mu)C_{i}$ and $B_{i}<r_{i}D_{i}$.}

It is straightforward to verify that Lemma \ref{lemma:MixedMPE2}
continues to hold under this more general model, because the only
properties of the payoff-relevant parameters we used in the proof
of this lemma are that (a) $w(x)>l_{i}$ for all $x$, and (b) $w(x)>\mathbb{E}\left[\int_{0}^{t}e^{-rs}\pi(X_{s})ds+\right.$
$\left.e^{-rt}w(X_{t})\,|X_{0}=x\right]$ for all $x\in\mathcal{X}$
and $t$. It thus follows that if Condition \ref{Ass-Str-Stability}
is satisfied and the parameters $\left\{ r_{i},\pi_{i}(\cdot),w_{i}(\cdot),l_{i}(\cdot)\right\} _{i\in\{1,2\}}$
are such that the thresholds $\theta_{1}^{*}\neq\theta_{2}^{*}$,
then the game admits no mixed-strategy MPE.

\section{Online Appendix: Singular Strategy \label{sec:Singular}}

In this section, we provide justification for precluding a singularly
continuous component when a strategy is expressed as a cumulative
distribution function of exit time.

We first note that firm $i$'s strategy can be alternatively expressed
as a non-decreasing and right-continuous process $A_{i}=\{A_{i,t}\}_{t\geq0}$
that ranges in the interval $[0,\infty]$; it can be transformed into
a cumulative distribution function $G_{i}(t)$ of firm $i$'s exit
timing by letting $\int_{0}^{t}\frac{dG_{i}(s)}{1-G_{i}(s)}=A_{i,t}$.
Such a process can be decomposed into three components as follows:
\[
A_{i,t}=\int_{0}^{t}\lambda_{i,s}ds+\int_{0}^{t}dL_{i,s}+\sum_{0<u\le t}\Delta A_{i,u}\:,
\]
where the first component $\int_{0}^{t}\lambda_{i,s}ds$ is the absolutely
continuous (in time) part, the second component $\int_{0}^{t}dL_{i,s}$
is the singularly continuous (in time) part, and the third component
$\sum_{0<u\le t}\Delta A_{i,u}$ is the discontinuous part.

We now argue that the singularly continuous part of $A_{i,t}$ is
absent, i.e., $\int_{0}^{t}dL_{i,s}=0$, in a mixed strategy MPE.
Suppose on the contrary that there exists a point $y\in\mathcal{X}\subseteq\mathbb{R}$
at which the increase $dA_{i,t}$ is singularly continuous. Let $L_{i}(t)$
denote the singularly continuous component of $A_{i}$. Then, from
the definition provided by \citet{KS1991}, we can express it as
\[
L_{i}(t)=\lim_{\epsilon\downarrow0}\int_{0}^{t}\frac{f(X_{s})}{\epsilon}\mathbf{1}_{\{y-\epsilon<X_{s}<y+\epsilon\}}ds\;.
\]
For example, $L_{i}(t)$ reduces to the local time if $f(\cdot)=1/2$
and if $X$ is a Wiener process. Note that we can think of $L_{i}(\cdot)$
as resulting from a very large exit rate $f(X_{t})/\epsilon$ at and
around $y$. Also, for any $\epsilon>0$, we can think of $(y-\epsilon,y+\epsilon)$
as a mixed strategy exit region with an exit rate of $f(X_{t})/\epsilon$.
Then because firm $i$ must be indifferent between exit at time $0$
and exit at an infinitesimal time $dt$, firm $-i$ must also have
mixed strategy exit region in $(y-\epsilon,y+\epsilon)$ with 
\[
\lambda_{-i}(x)=\frac{rl_{i}-\pi(x)}{w(x)-l_{i}}\:.
\]
Similarly, by a symmetric argument, firm $i$'s exit rate should be
\[
\lambda_{i}(x)=\frac{rl_{-i}-\pi(x)}{w(x)-l_{-i}}\:,
\]
which cannot be arbitrarily large as $f(X_{t})/\epsilon$ as $\epsilon\downarrow0$.
Therefore, such a singularly continuous component cannot exist in
a mixed strategy MPE.

\section{Online Appendix: Proof of Lemma \ref{lemm:Opt-Stop} \label{sec:Lemma1}}

\noindent The proof of this lemma is available in \citet{Alvarez2001},
but here, we provide a sketch of the proof based on the verification
theorem \citep[Theorem 10.4.1]{Oksendal03}. To that end, we will
use the optimality conditions, which are known as ``value matching''
and ``smooth pasting'' conditions \citep{Samuelson1965,McKean1965,Merton1973}.

First, the state space $\mathcal{X}$ must be the union of $C:=\{x\in\mathcal{X}:V_{i}^{*}(x)>l_{i}\}$
and $\Gamma:=\{x\in\mathcal{X}:V_{i}^{*}(x)=l_{i}\}$, which are mutually
exclusive: This is because (1) $X$ is a time-homogeneous process
and the time horizon is infinite, and (2) the value function $V_{i}^{*}(\cdot)$
from an optimal stopping policy must be always no less than the reward
$l_{i}$ from stopping immediately. Hence, the problem to find an
optimal stopping policy can be reduced to identify $C$ or $\Gamma$.

Next, we find the differential equation that $V_{i}^{*}(x)$ must
satisfy \emph{if} $x\in C$. Note that the optimal value function
$V_{i}^{*}(x)$ is the maximum of the reward from waiting an instant
and the reward from stopping immediately. For any $x\in C$, therefore,
the optimal stopping policy is to wait an instant $dt$, and hence,
the optimal value function must satisfy the following equation:
\begin{align}
V_{i}^{*}(x) & =\pi(x)dt+(1-rdt)\mathbb{E}^{x}[V_{i}^{*}(x)+dV_{i}^{*}(X_{t})]\:.\label{eq:Pf-Conti}
\end{align}
Then applying Ito formula to $V_{i}^{*}(X_{t})$ and using $\mathbb{E}^{x}[dB_{t}]=0$
yields
\begin{align}
\mathbb{E}^{x}[dV_{i}^{*}(X_{t})] & =[\mu(x)V_{i}^{*'}(x)+\frac{1}{2}\sigma^{2}(x)V_{i}^{*''}(x)]dt\:.\label{eq:Pf-dVi}
\end{align}
By plugging \eqref{eq:Pf-dVi} into \eqref{eq:Pf-Conti} and ignoring
the term smaller than $dt$, we have
\begin{align*}
V_{i}^{*}(x) & =\pi(x)dt+V_{i}^{*}(x)+[-rV_{i}^{*}(x)+\mu(x)V_{i}^{*'}(x)+\frac{1}{2}\sigma^{2}(x)V_{i}^{*''}(x)]dt\:,
\end{align*}
from which we obtain the following second-order linear differential
equation:
\begin{align}
\frac{1}{2}\sigma^{2}(x)V_{i}^{*''}(x)+\mu(x)V_{i}^{*'}(x)-rV_{i}^{*}(x) & =-\pi(x)\:.\label{eq:Pf-VarEq}
\end{align}
Thus, $V_{i}^{*}(\cdot)$ can be obtained by solving the differential
equation \eqref{eq:Pf-VarEq}. In fact, it can be seen from a series
of algebra with the relation \eqref{eq:Pf-LemmaA1-R} that the function
$R(\cdot)+A\phi(\cdot)$ with some constant $A\in\mathbb{R}$ is a
solution to \eqref{eq:Pf-VarEq}, and hence, we can guess $V_{i}^{*}(x)=R(x)+A\phi(x)$
with some constant $A$. 

Intuitively, firm $i$ must find it optimal to exit and receive his
outside option $l_{i}$ as soon as the state $X$ hits some lower
threshold $\theta_{i}$. Hence, assume at the moment that the optimal
stopping policy is given as $\tau^{*}:=\inf\{t\geq0:X_{t}^{x}\leq\theta_{i}\}$,
which implies that $\theta_{i}$ is the boundary point of the region
$C$. Now, we state the value matching condition and the smooth pasting
condition, which results in two boundary conditions to the boundary
value problem \eqref{eq:Pf-VarEq} with the free boundary $\theta_{i}$:
\begin{align}
V_{i}^{*}(\theta_{i}) & =R(\theta_{i})+A\phi(\theta_{i})=l_{i}\label{eq:Pf-ValMatch}\\
V_{i}^{*'}(\theta_{i}) & =R'(\theta_{i})+A\phi'(\theta_{i})=0\:.\label{eq:Pf-SmthPaste}
\end{align}

The value matching condition \eqref{eq:Pf-ValMatch} and the smooth
pasting condition \eqref{eq:Pf-SmthPaste} are the conditions that
$V_{i}^{*}(\cdot)$ must satisfy at the boundary $\theta_{i}$ of
$C$. We can first obtain $A=[l_{i}-R(\theta_{i})]/\phi(\theta_{i})=\beta_{i}(\theta_{i})$
from \eqref{eq:Pf-ValMatch}. Then the condition \eqref{eq:Pf-SmthPaste}
is equivalent to
\begin{align*}
0 & =R^{'}(\theta_{i})+\frac{l_{i}-R(\theta_{i})}{\phi(\theta_{i})}\phi{}^{'}(\theta_{i})\\
 & =\frac{R^{'}(\theta_{i})\phi(\theta_{i})+[l_{i}-R(\theta_{i})]\phi{}^{'}(\theta_{i})}{\phi(\theta_{i})}=-\phi(\theta_{i})\beta_{i}^{'}(\theta_{i})\:.
\end{align*}
Because $\phi(\cdot)>0$, it can be seen from Lemma \ref{lemm:Aux_1}
that this condition is satisfied if and only if $\theta_{i}=\theta_{i}^{*}$,
which implies that $A=\beta_{i}(\theta_{i}^{*})$.

Lastly, it can be easily verified that $R(x)+\beta_{i}(\theta_{i}^{*})\phi(x)\geq l_{i}$
for $\forall x\geq\theta_{i}^{*}$ and $\pi(x)<rl_{i}$ for $\forall x\leq\theta_{i}^{*}<x_{ci}$.
By the verification theorem \citep[Theorem 10.4.1]{Oksendal03}, therefore,
the proposed value function $R(\cdot)+\beta_{i}(\theta_{i}^{*})\phi(\cdot)$
is, in fact, the optimal value function $V_{i}^{*}(\cdot)$, as desired.
\eproof

\section{Online Appendix: Extension to (Non-Markov) Subgame Perfect Equilibria
\label{sec:Online-nonMarkov}}

In this appendix, we allow firms to condition their decision at $t$
on the entire history $h^{t}=\{X_{s}\}_{s\leq t}$ (as opposed to
only the current state, $X_{t}$). We show that if the firms have
heterogeneous exit payoffs (\emph{i.e.,} $l_{1}<l_{2}$), then subject
to a set of restrictions on their strategies, the game admits no mixed-strategy
Subgame Perfect equilibrium.

We first extend the strategy to accommodate the history dependence.
At every moment, given the history $h^{t}=\left\{ X_{s}\right\} _{s\leq t}$
and conditional on the game not having ended, each firm chooses (probabilistically)
whether to exit. Formally, each firm $i$ chooses
\begin{enumerate}
\item [i.]a set of histories $I_{i}$ (or an \emph{exit region}) such that
it exits with probability $1$ if $h^{t}\in I_{i}$,
\item [ii.]a set of stopping time and exit probability pairs, denoted by
$\mathcal{P}_{i}=\{(\tau_{i,n},p_{i,n})\}_{n=1}^{\infty}$, such that
it exits at $t=\tau_{i,n}$ with probability $p_{i,n}\in(0,1)$, and
\item [iii.]a non-negative process $\Lambda_{i}=\{\lambda_{i,t}\}_{t>0}$,
which represents the firm's hazard rate of exit at $t$.
\end{enumerate}
We assume that $\tau_{I}:=\inf\{t:h^{t}\in I_{i}\}$ is a stopping
time. We also assume that each $p_{i,n}$ and the process $\lambda_{i,t}$
is progressively measurable with respect to $\mathcal{F}_{\tau_{i,n}}$
and $\mathcal{F}_{t}$, respectively, and so the exit probability
$p_{i,n}$ at $\tau_{i,n}$, and the exit rate $\lambda_{i,t}$ may
depend on the entire history of $X$ up to $\tau_{i,n}$ and $t$,
respectively. Then we can represent firm $i$'s strategy as the
three-tuple $a_{i}=(I_{i},\mathcal{P}_{i},\Lambda_{i})$, and $\left\{ a_{1},a_{2}\right\} $
is a \emph{strategy profile}. As each firm's decision at $t$ can
be conditioned on the entire history up to $t$, it is without loss
of generality to assume that each firm chooses its strategy at time
$0$. Intuitively, during any \emph{small} interval $[t,t+dt)$, firm
$i$ exits with probability
\[
\rho_{i,t}=\begin{cases}
1 & \text{ if }h^{t}\in I_{i}\text{ ,}\\
p_{i,n}+(1-p_{i,n})\lambda_{i,t}dt & \text{ if }t=\tau_{i,n}\text{ , and}\\
\lambda_{i,t}dt & \text{ otherwise.}
\end{cases}
\]

If firm $i$ does not exit with probability $1$ after any history,
we write $I_{i}=\emptyset$. If it does not choose any stopping time
-- exit probability pairs, we write $\mathcal{P}_{i}=\emptyset$.
If it does not exit with a positive hazard rate (\emph{i.e.,} $\lambda_{i,t}=0$
for all $t$ almost surely, hereafter a.s), we write $\Lambda_{i}=\mathbf{0}$.
Finally, we say that firm $i$'s strategy is\emph{ pure} if $\mathcal{P}_{i}=\emptyset$
and $\Lambda_{i}\equiv\mathbf{0}$, and it is \emph{mixed} otherwise.

Next, we write each firm's payoff as a function of an arbitrary strategy
profile. Fix a strategy profile $\left\{ a_{1},a_{2}\right\} $ and
history $h^{t}$, and define $\tau_{i,0}:=\inf\left\{ s\ge0:h^{s}\in I_{i}\right\} $
and $p_{i,0}:=1$. The survival probability, that is, the probability
that the game does not end during $[t,u)$ is given by
\[
S_{t,u}:=e^{-\int_{t}^{u}(\lambda_{1,s}+\lambda_{2,s})ds}\prod_{\{n,m\ge0:t\le\tau_{1,n}<u,t\le\tau_{2,m}<u\}}(1-p_{1,n})(1-p_{2,m})\:.
\]
Firm $i$'s payoff at time $t$ (conditional on the game not having
ended) can be written as
\begin{align}
V_{i}(h^{t};a_{1},a_{2})= & \mathbb{E}\left[\int_{t}^{\infty}e^{-r(s-t)}S_{t,s}[\pi(X_{s})+\lambda_{i,s}l_{i}+\lambda_{j,s}w(X_{s})]ds\right.\nonumber \\
 & +\sum_{n\geq0}S_{t,\tau_{i,n}}e^{-r(\tau_{i,n}-t)}p_{i,n}l_{i}+\sum_{m\geq0}S_{t,\tau_{j,m}}e^{-r(\tau_{i,m}-t)}p_{j,m}w(X_{\tau_{j,m}})\nonumber \\
 & \left.-\frac{1}{2}\sum_{n,m\geq0}\mathbb{I}_{\left\{ \tau_{i,n}=\tau_{j,m}\right\} }S_{t,\tau_{i,n}}e^{-r(\tau_{i,n}-t)}p_{i,n}p_{j,m}\left(l_{i}+w(X_{\tau_{j,m}})\right)\right]\,.\label{eq:Payoff-2}
\end{align}
The first line represents the firm's discounted flow payoff with survival
chances taken into account, plus the reward from the end of the game
through the exit rate by either firm. The second line captures the
payoff from either firm's instantaneous exit probability, while the
third line accounts for the possibility of simultaneous exit and the
double counting from the second line. The dependence of the strategies
$a_{i}$ on history $h^{t}$ is muted for expositional simplicity.

A strategy profile $\left\{ a_{1}^{*},a_{2}^{*}\right\} $ is a Subgame
Perfect equilibrium (hereafter SPE) if 
\[
V_{i}(h^{t};a_{i}^{*},a_{-i}^{*})\ge V_{i}(h^{t};a_{i},a_{-i}^{*})
\]
for each firm $i$, every history $h^{t}$, and every strategy $a_{i}$.\footnote{\citet{DR1993} uses a similar formulation in a class of stopping
time games.}

Recall that firm $i$'s strategy can be summarized by the three-tuple
$(I_{i},\mathcal{P}_{i},\Lambda_{i})$, where $I_{i}$ is a set of
histories such that firm $i$ exits instantaneously whenever $h^{t}\in I_{i}$,
$\mathcal{P}_{i}=\{\tau_{i,n},p_{i,n}\}_{n=1}^{\infty}$ is a collection
of stopping times and corresponding exit probabilities, and $\Lambda_{i}=\left\{ \lambda_{i,t}\right\} _{t\geq0}$
is a non-negative process. To make it explicit that $\lambda_{i,t}$
can depend on the entire history $h^{t}$, we will sometimes write
$\lambda_{i}(h)$ to denote firm $i$'s exit rate when $h^{t}=h$.\medskip{}

To help the reader visualize an SPE with history-dependent mixed strategies,
we present an example when the firms are homogeneous. \begin{example}\label{ex:mixed SPE}Suppose
that $l_{1}=l_{2}$ (and so $\theta_{1}^{*}=\theta_{2}^{*}$ by Lemma
\ref{lemm:Opt-Stop}). Fix any $q\in(0,1)$, and consider the strategies
$a_{1}=\left(\emptyset,\left\{ \tau_{1},q\right\} ,\left\{ \lambda_{t}\right\} _{t\geq0}\right)$
and $a_{2}=\left(\emptyset,\emptyset,\left\{ \lambda_{t}\right\} _{t\geq0}\right)$,
where $\tau_{1}=\inf\left\{ t\geq0\,:\,X_{t}\leq\theta_{1}^{*}\right\} $,
and 
\[
\lambda_{t}:=\mathbb{I}_{\left\{ X_{t}\leq\theta_{1}^{*}\right\} }\frac{rl_{1}-\pi(X_{t})}{w(X_{t})-l_{1}}\,.
\]

\noindent Then $\left\{ a_{1},a_{2}\right\} $ constitutes a (non-Markov)
mixed strategy SPE.\end{example}

In this example, both firms remain in the market until the first time
that $X_{t}\in(\alpha,\theta_{1}^{*}]$. At that moment, firm 1 exits
with instantaneous probability $q$. From that time onwards, whenever
$X_{t}\leq\theta_{1}^{*}$, each firm exits with rate $\lambda_{t}$,
which is chosen to make its opponent indifferent between remaining
in the market and exiting. This strategy profile is non-Markov because
firm $1$ exits with probability $q$ only at the first time that
$X_{t}$ hits $\theta_{1}^{*}$. Indeed, this is the stochastic analog
of the mixed-strategy equilibrium characterized in \citet{Hendricks1988}
when the game is deterministic (\emph{i.e.,} $\sigma(\cdot)\equiv0$),
and in \citet{Tirole1988} when $\mu(\cdot)\equiv\sigma(\cdot)\equiv0$
and $X_{0}$ satisfies $\pi(X_{0})<rl_{i}$.\footnote{The reader is referred to \citet{Steg2015} for a proof that the proposed
strategies indeed constitute an SPE.}

We note that the strategy profile $\left\{ a_{1},a_{2}\right\} $
in Example \ref{ex:mixed SPE} cannot constitute a mixed strategy
SPE if the firms are heterogeneous, i.e., $l_{1}<l_{2}$. This is
because if firm $1$ does not exit at the first time $X_{t}\leq\theta_{1}^{*}$,
which happens with probability $1-q>0$, then both firms' strategies
are Markov afterwards. In the subgame firm $1$ stays at the first
time $X_{t}\leq\theta_{1}^{*}$, therefore, the strategy profile $\left\{ a_{1},a_{2}\right\} $
constitutes a mixed strategy MPE, which has been proved impossible
in Section \ref{sec:MPE}; in the rest of this section, we will formalize
and generalize this argument. In addition, we argue that simple variations
of the strategy profile $\left\{ a_{1},a_{2}\right\} $ in Example
\ref{ex:mixed SPE} would not constitute a mixed strategy SPE if $l_{1}<l_{2}$.
For instance, one may consider a strategy where firm $1$ exits with
probability $q$ whenever $X_{t}$ hits $\theta_{1}^{*}$ only from
above. This strategy, however, is indistinguishable from a strategy
where firm $1$ exits whenever $X_{t}$ hits $\theta_{1}^{*}$ from
either side because $X$ is a regular diffusion process. As another
variation, one may consider the adjustment of exit probability $q$
depending on the history $h^{t}$ while still keeping the threshold
rule: Whenever $X_{t}$ hits $\theta_{1}^{*}$ from above, firm $1$
exits with probability $q(h^{t})$ depending on the history $h^{t}$.
However, no matter how carefully $q(h^{t})>0$ is chosen, firm $2$
would not exit with the rate $\lambda_{t}$ near but below $\theta_{1}^{*}$
because of the chance to become the winner soon.

\medskip{}

We impose three restrictions on the firms' strategies, all of which
are satisfied by the strategies in the mixed-strategy equilibria that
appear in the extant literature; see for instance, \citet{Tirole1988},
\citet{Hendricks1988}, \citet{Levin2004}, \citet{Steg2015}, and
Example 1 above. The first is that the exit regions $I_{1}=I_{2}=\emptyset$;
\emph{i.e.}, neither firm exits with probability one following any
history. The second is that $\mathcal{P}_{i}=\{\tau_{i,n},p_{i,n}\}_{n=1}^{N_{i}}$
for some $N_{i}<\infty$ and $\tau_{i,n}<\infty$ a.s for all $i$
and $n$. That is, the number of events of instantaneous exit is finite
and the stopping times are finite a.s. We discuss the role of these
assumptions before Lemma \ref{lemma:common-support-SPE}. Finally,
we define
\[
H_{i}=\left\{ h\,:\,\lambda_{i}(h)>0\right\} 
\]
and impose the regularity condition that $H_{i}$ is an open set with
respect to the following metric $d(\cdot,\cdot)$ on the space of
histories: 
\begin{equation}
d(h_{1},h_{2}):=\sqrt{\max\left\{ \max_{s\in[0,t]}\left|X_{s}^{1}-X_{s}^{2}\right|,\sup_{s\in(t,t^{\prime}]}\left|X_{t}^{1}-X_{s}^{2}\right|\right\} +(t-t^{\prime})^{2}}\,,\label{eq:metric}
\end{equation}
where $h_{1}^{t}$ and $h_{2}^{t'}$ are two histories with $t\le t'$.
This is a generalization of the uniform metric, which takes into account
the different lengths of the histories. It is straightforward to verify
that $d(\cdot,\cdot)$ fits the definition of a metric, and hence,
we can define open sets of histories with respect to this metric.
Note that $H_{i}$ is the counterpart of $\Gamma_{i}$ defined in
(\ref{eq:supp-lambda}) when strategies are non-Markov: It comprises
the histories in which firm $i$ randomizes between remaining in the
market and exiting. With the assumption that $H_{i}$ is an open set,
firm $i$'s strategy is smooth in the sense that if $\lambda(h^{t})>0$,
then for any other history $h^{t'}$ close enough to $h^{t}$, we
have $\lambda_{i}(h^{t'})>0$. We summarize these conditions below.\begin{condition}\label{Ass-Restrictions-SPE}
Assume that each firm $i$'s strategy $a_{i}=\left(I_{i},\{\tau_{i,n},p_{i,n}\}_{n=1}^{N_{i}},\lambda_{i}(h)\right)$,
where

\noindent (i) $I_{i}=\emptyset$ for each $i$,

\noindent (ii) $N_{i}<\infty$ and $\tau_{i,n}<\infty$ a.s for all
$i$ and $n$, and

\noindent (iii) $H_{i}=\left\{ h\,:\,\lambda_{i}(h)>0\right\} $ is
an open set. \end{condition}

\noindent Conditions \ref{Ass-Restrictions-SPE}(i) and (ii) ensure
that there exists a stopping time $\overline{\tau}:=\max_{n}\{\tau_{1,n},\tau_{2,n}\}$
such that exit after this time occurs only via the hazard rate, which
implies that there exists a subgame (starting at $\overline{\tau}$)
that is reached with positive probability, in which exit occurs only
via the hazard rates.\footnote{Put differently, each firm exits during any interval of length $dt$
with probability $\lambda_{i}(h^{t})dt$ for all $t>\overline{\tau}$.
The restriction that $\tau_{i,n}<\infty$ a.s simplifies the exposition
by ensuring that $\max_{n}\{\tau_{1,n},\tau_{2,n}\}<\infty$ a.s.
It can be relaxed, and a proof of the results in this section absent
this condition is available upon request.}\footnote{Note that, with Condition \ref{Ass-Restrictions-SPE}(i), we can focus
on ``proper'' mixed strategies, in which the probability of exit
is always less than one, and this is the case in all the known examples
of mixed-strategy equilibria in the literature.} Let $H_{i}^{\overline{\tau}}$ denote the set of histories in that
subgame (\emph{i.e.,} the histories $h^{t}\supseteq h^{\overline{\tau}}$)
such that $\lambda_{i}(h)>0$. The following lemma shows that in a
(mixed-strategy) SPE, the firms must randomize between remaining in
the market and exiting over a common set of histories in that subgame.\begin{lemma}\label{lemma:common-support-SPE}
Suppose Condition \ref{Ass-Restrictions-SPE} is satisfied and $\left\{ a_{1},a_{2}\right\} $
constitutes a mixed-strategy SPE. Then $H_{1}^{\overline{\tau}}=H_{2}^{\overline{\tau}}$.
\end{lemma}

\noindent \textbf{Proof of Lemma} \textbf{\ref{lemma:common-support-SPE}}:
Define the stopping time $\overline{\tau}:=\max_{n}\left\{ \tau_{1,n},\tau_{2,n}\right\} $,
and note that it is finite a.s by Condition \ref{Ass-Restrictions-SPE}(ii).
Let $H_{i}^{\overline{\tau}}:=\{h\supseteq h^{\overline{\tau}}:\lambda_{i}(h)>0\}$
denote the set of histories $h$ in the subgame $h^{\overline{\tau}}$
in which firm $i$ exits with hazard rate $\lambda_{i}(h)>0$. Next,
given a history $h^{v}$ with $v>\overline{\tau}$, define $H_{i}^{v}:=\{h^{t}\supseteq h^{v}:\lambda_{i}(h^{t})>0\}$.
Then because $H_{i}^{\overline{\tau}}=\cup_{v>\overline{\tau}}H_{i}^{v}$,
it is enough to show that $H_{1}^{v}=H_{2}^{v}$ for any given history
$h^{v}$ with $v>\overline{\tau}$. Now, towards a contradiction,
fix some history $h^{v}$ with $v>\overline{\tau}$ and suppose that
there exists an open set $C_{1}\subset H_{1}^{v}$ such that $C_{1}\cap H_{2}^{v}=\emptyset$.
Pick a history $h^{t}\in C_{1}$ and define $\tau(\delta):=\inf\{s\geq t:X_{s}\notin(X_{t}-\delta,X_{t}+\delta)\}$,
i.e., $\tau(\delta)$ is the exit time of a $\delta$-neighborhood
of $X_{t}$.

\begin{onehalfspace}
We first note that because $C_{1}$ is open and $h^{t}\in C_{1}$,
there exists some $\epsilon>0$ such that any history $h$ with $d(h^{t},h)<\epsilon$
must belong to $C_{1}$ (i.e., a small enough neighborhood of $h^{t}$
must be contained in $C_{1}$). Next, based on the definition of the
metric $d(\cdot,\cdot)$ on the space of histories of $X$, there
must exist some $\delta>0$ and some time $u>t$ such that $d(h^{t},h^{\tau(\delta)\land u})<\epsilon$,
which implies that $h^{\tau(\delta)\land u}\in C_{1}$ by our first
note. This is because any history $h^{t'}$in the subgame $h^{t}$
(i.e., a history $h^{t'}\supseteq h^{t}$) must be close enough to
$h^{t}$ with respect to the metric $d(\cdot,\cdot)$ as long as $X_{t'}$
(the value of $X$ in the history $h^{t'}$) is not too much different
from $X_{t}$ and $t'-t$ is not too large. Then for any $s\in[t,\tau(\delta)\land u]$,
we have $d(h^{t},h^{s})\leq d(h^{t},h^{\tau(\delta)\land u})<\epsilon$,
which implies that $h^{s}\in C_{1}$ for any $s\in[t,\tau(\delta)\land u]$.
Hence, for any $s\in[t,\tau(\delta)\land u]$, we have $h^{s}\notin H_{2}$
because $C_{1}\cap H_{2}=\emptyset$.

Next, recall that $h^{t}\in C_{1}\subseteq H_{1}^{v}$ and $v>\overline{\tau}$,
which implies that $s>\overline{\tau}$ for any $s\in[t,\tau(\delta)\land u]$.
Let $H^{\tau(\delta)\land u}:=\{h:h\supseteq h^{\tau(\delta)\land u}\}$
and $H^{t}:=\{h:h\supseteq h^{t}\}$ be the sets of all the histories
that contain $h^{\tau(\delta)\land u}$ and $h^{t}$ respectively.
Then because $h^{s}\notin H_{2}^{v}$ and $s>\overline{\tau}$ for
all $s\in[t,\tau(\delta)\land u]$, firm $2$ does not exit at all
within this time interval, which means that
\[
V_{1}(h^{t};\left(H^{\tau(\delta)\land u},\mathcal{P}_{1},\mathbf{0}\right),a_{2})-V_{1}(h^{t};\left(H^{t},\mathcal{P}_{1},\mathbf{0}\right),a_{2})=\mathbb{E}\biggl[\int_{t}^{\tau(\delta)\land u}\left(\pi(X_{s})-rl_{1}\right)e^{-r(s-t)}ds\biggr]<0\:,
\]
where the inequality holds because $X_{s}<\theta_{i}^{*}$ is a necessary
condition for a mixed strategy region for firm $1$, in which case
$\pi(X_{s})-rl_{1}<0$. Here $\mathbf{0}$ means that $\lambda_{i}(h)=0$
for any history $h$. This, however, contradicts that $h^{\tau(\delta)\land u}\in C_{1}$
because the payoff from exiting at time $\tau(\delta)\land u$ is
strictly less than that from exiting at time $t$, which implies that
$\lambda_{1}(h^{\tau(\delta)\land u})>0$ is not a best response to
$a_{2}$.\eproof
\end{onehalfspace}

This is a counterpart of Lemma \ref{lemma:MixedMPE2} when strategies
are not constrained to be Markov. It is helpful to convey the intuition
with a heuristic derivation. Towards a contradiction, suppose that
there exists a non-empty open set of histories $C\subset H_{i}^{\overline{\tau}}\backslash H_{-i}^{\overline{\tau}}$.
Pick a history $h^{t}\in C$. Then $V_{i}(h^{\prime};a)=l_{i}$ for
any $h^{\prime}$ in the vicinity of $h^{t}$; this is because any
history $h^{\prime}$ close enough to $h^{t}$ must also belong to
$C\subset H_{i}^{\overline{\tau}}$, and if so, we must have $\lambda_{i}(h^{\prime})>0$,
which requires $V_{i}(h^{\prime};a)=l_{i}$ for any mixed-strategy
equilibrium $a=(a_{i},a_{-i})$. Hence, there is an infinitesimal
time $\Delta t$ such that firm $i$ will be indifferent between immediate
exit at $t$ and exit at $t+\Delta t$:
\[
l_{i}=\pi(X_{t})\Delta t+\exp(-r\Delta t)l_{i}+o(\Delta t)=l_{i}+[\pi(X_{t})-rl_{i}]\Delta t+o(\Delta t)\,,
\]
where we have used that an exit after an infinitesimal time earns
$l_{i}$ because $h^{t+\Delta t}\in C$. Then the indifference equation
leads to $\pi(X_{t})=rl_{i}$, which contradicts the fact that $\pi(X_{t})<rl_{i}$
for any $h^{t}\in H_{i}$. Therefore, we can conclude that $C$ is
empty, and so $H_{1}^{\overline{\tau}}=H_{2}^{\overline{\tau}}$.

\medskip{}

The following theorem shows that if the firms have heterogeneous outside
options, then subject to Condition \ref{Ass-Restrictions-SPE}, the
game admits no mixed-strategy SPE; \emph{i.e.,} there exists no SPE
such that $H_{1}\cup H_{2}\neq\emptyset$.\begin{theorem}\label{Thm-SPE}Suppose
that each firm's strategy must satisfy Condition \ref{Ass-Restrictions-SPE}.
If $l_{1}<l_{2}$, then no mixed-strategy SPE\emph{ }exists\emph{.}\end{theorem}

\noindent \textbf{Proof of Theorem \ref{Thm-SPE}.} Towards a contradiction,
suppose that $l_{1}<l_{2}$, the strategies $a_{1}$ and $a_{2}$
satisfy Condition \ref{Ass-Restrictions-SPE}, and the strategy profile
$\{a_{1},a_{2}\}$ constitutes a mixed-strategy SPE. By Lemma \ref{lemma:common-support-SPE},
we have $H_{1}^{\overline{\tau}}=H_{2}^{\overline{\tau}}$, so define
$H^{\overline{\tau}}:=H_{1}^{\overline{\tau}}$.

Fix a finite time $t$ such that $t>\overline{\tau}$ a.s and $X_{t}\geq\max\{\theta_{1}^{*},\theta_{2}^{*}\}$
(reached with positive probability), and let $\tau_{H}:=\inf\{s>t:h^{s}\in H^{\overline{\tau}}\}$
denote the first hitting time of the mixed-strategy region after time
$t$. We first note that $\lambda_{i}$ is defined as a progressively
measurable non-negative process (See Section \ref{subsec:Strategies}),
i.e., $\lambda_{i}=\{\lambda_{i,s}\}_{s\geq0}$ where $\lambda_{i,s}$
is firm $i$'s hazard rate of exit at time $s$. Then based on the
definitions of $H^{\overline{\tau}}$ and $\tau_{H}$, we have
\[
\tau_{H}=\inf\{s>t:\lambda_{i,s}>0\}\,,
\]
which is a stopping time with respect to $\mathcal{F}_{X}$.

Next, we have already established that $H$ does not intersect with
the region in which $X_{t}>\theta_{i}^{*}$ for either $i$. Because
$X_{t}\geq\max\{\theta_{1}^{*},\theta_{2}^{*}\}$, we have $\tau_{H}\ge\tau_{i}^{*}:=\inf\{s>t:X_{s}\leq\theta_{i}^{*}\}$
for each $i$ a.s, that is, the first hitting time of $H^{\overline{\tau}}$
is at least as long as the hitting time of $(a,\theta_{1}^{*}]$,
as well as $(a,\theta_{2}^{*}]$.

We now show that $\tau_{H}=\tau_{i}^{*}$ a.s for both $i$. Towards
a contradiction, suppose that $\tau_{H}\neq\tau_{i}^{*}$ with positive
probability. At $\tau_{H}$, firm $i$'s payoff 
\[
V_{i}(h^{\tau_{H}};a_{i},a_{-i})=l_{i}=V_{i}(h^{\tau_{H}};\left(H^{\overline{\tau}},\mathcal{P}_{i},\mathbf{0}\right),a_{-i})=V_{i}(h^{\tau_{H}};\left(H^{\overline{\tau}},\mathcal{P}_{i},\mathbf{0}\right),\left(\emptyset,\mathcal{P}_{-i},\mathbf{0}\right))\,.
\]
The first equality follows from the fact that firm $i$'s payoff must
be equal to $l_{i}$ at $\tau_{H}$ by definition of a mixed-strategy
equilibrium. The third term represents firm $i$'s payoff function
if it exits with probability one whenever $h^{t}\in H^{\overline{\tau}}$
and its opponent's strategy is $a_{-i}$, whereas the last term represents
firm $i$'s payoff function if it exits with probability one whenever
$h^{t}\in H^{\overline{\tau}}$ and its opponent never exits. The
third equality follows because $\tau_{H}>\overline{\tau}$ (so $\mathcal{P}_{-i}$
has no impact on firm $i$'s payoff) and $\lambda_{-i,\tau_{H}}=0$
by the definition of $\tau_{H}$ and $\overline{\tau}$. The second
and third equalities assert that firm $i$'s payoff at $\tau_{H}$
under the specified strategy profiles equals $l_{i}$.

Therefore, if $\tau_{H}\neq\tau_{i}^{*}$ with non-zero probability,
then there exists a history $h^{t}$ (where $t$ has been fixed in
the beginning of the proof) such that 
\begin{align*}
V_{i}(h^{t};\left(H^{\overline{\tau}},\mathcal{P}_{i},\mathbf{0}\right),a_{-i}) & =V_{i}(h^{t};\left(H^{\overline{\tau}},\mathcal{P}_{i},\mathbf{0}\right),\left(\emptyset,\mathcal{P}_{-i},\mathbf{0}\right))\\
 & <V_{i}(h^{t};\left(H^{\tau_{i}^{*}},\mathcal{P}_{i},\mathbf{0}\right),\left(\emptyset,\mathcal{P}_{-i},\mathbf{0}\right))=V_{i}(h^{t};\left(H^{\tau_{i}^{*}},\mathcal{P}_{i},\mathbf{0}\right),a_{-i})\,,
\end{align*}
where $H^{\tau_{i}^{*}}:=\{h^{t}:X_{t}\leq\theta_{i}^{*}\}$ is the
set of histories in which $X_{t}\leq\theta_{i}^{*}$. Here the inequality
follows because by Lemma \ref{lemm:Opt-Stop}, exiting whenever $X_{t}\leq\theta_{i}^{*}$
is firm $i$'s unique best response to an opponent who never exits,
and $\tau_{H}\neq\tau_{i}^{*}$ with positive probability. The first
equality follows because $t>\overline{\tau}$ (so $\mathcal{P}_{-i}$
has no impact on firm $i$'s payoff) and $\lambda_{-i,\tau_{H}}=0$.
The second equality follows because $t>\overline{\tau}$ and $\tau_{H}\ge\tau_{i}^{*}$
a.s. This contradicts the assumption that $(a_{1},a_{2})$ is an SPE,
and hence we conclude that $\tau_{H}=\tau_{i}^{*}$ a.s. However,
this is not possible if $\theta_{1}^{*}\neq\theta_{2}^{*}$, or equivalently,
if $l_{1}<l_{2}$, yielding a contradiction.

Finally, note that any such $t$ is reached with positive probability
by Condition \ref{Ass-Restrictions-SPE}(i) and because $X$ is irreducible.
Therefore, it follows that no mixed-strategy SPE exists. \eproof

We briefly give a summary of the proof here. First recall that $H_{1}^{\overline{\tau}}=H_{2}^{\overline{\tau}}$
by Lemma \ref{lemma:common-support-SPE}. For any history in $H_{i}^{\overline{\tau}}$,
an indifference condition similar to (\ref{eq:Indifference}) must
be satisfied, which implies that $\lambda_{1,t}$ and $\lambda_{2,t}$
depend only on the current state, $X_{t}$. It is shown, using a similar
argument as in Lemma \ref{lemma:MixedMPE2}, that $H_{i}^{\overline{\tau}}$
must consist of histories such that $X_{t}\leq\theta_{i}^{*}.$ But
this implies that $H_{1}^{\overline{\tau}}\neq H_{2}^{\overline{\tau}}$
whenever $l_{1}<l_{2}$, a contradiction. Finally, because $H_{i}^{\overline{\tau}}$
is reached with positive probability (by Conditions \ref{Ass-Restrictions-SPE}(i)
and (ii)), it follows that no mixed-strategy SPE exists.
\end{document}